\documentclass[11pt]{article}
\usepackage[dvips]{graphicx}
\usepackage{latexsym}
\usepackage{amssymb}

\newcommand{\C}{\mathbb{C}}
\newcommand{\CP}{\mathbb{CP}}

\newcommand{\R}{\mathbb{R}}
\newcommand{\bbS}{\mathbb{S}}
\newcommand{\tbbS}{\tilde{\mathbb{S}}}
\newcommand{\Z}{\mathbb{Z}}
\newcommand{\X}{\mathbb{X}}

\newcommand{\cK}{{\mathcal{K}}}
\newcommand{\Lie}{{\mathcal{L}}}
\newcommand{\cM}{{\mathcal{M}}}
\newcommand{\cN}{{\mathcal{N}}}
\newcommand{\cPT}{{\mathcal{PT}}}
\newcommand{\PT}{{\mathcal{PT}}}
\newcommand{\half}{\frac{1}{2}}
\renewcommand{\d}{\mathrm{d}}

\newcommand{\Proof}{\noindent {\bf Proof: }}
\newcommand{\koniec}{\hfill  $\Box $\medskip}
\def\be{\begin{equation}}
\def\ee{\end{equation}}

\def\Sm{\Sigma}

\def\Om{\Omega}
\def\Th{\Theta}

\def\om{\omega}

\def\p{\partial}

\newcommand{\hook}{{\setlength{\unitlength}{11pt}   
                   \begin{picture}(.833,.8)
                   \put(.15,.08){\line(1,0){.35}}
                   \put(.5,.08){\line(0,1){.5}}
                   \end{picture}}}
\def\a{\alpha}
\def\ll{\lambda}

\def\OO{{\cal O}}

\newtheorem{theo}{Theorem}[section] 
\newtheorem{prop}[theo]{Proposition}  
\newtheorem{lemma}[theo]{Lemma}
\newtheorem{defi}[theo]{Definition}

\begin{document}
\title{Twistor theory of hyper-K{\"a}hler metrics with hidden symmetries}

\author{Maciej Dunajski\\
Department of Applied Mathematics and Theoretical Physics, 
University of Cambridge\\
Wilberforce Road,
Cambridge,
CB3 OWA,
UK,\\
\and
Lionel J. Mason\\ 
Mathematical Institute,
University of  Oxford \\
24-29 St Giles, Oxford OX1 3LB, UK.
}
\date{} 
\maketitle
\abstract {We review the hierarchy
for the hyper-K\"ahler equations and define a notion of symmetry
for solutions of this hierarchy. A
four-dimensional hyper-K\"ahler metric admits a hidden symmetry if it
embeds into a hierarchy with a symmetry.  It is shown that a
hyper-K\"ahler metric admits a hidden symmetry if it admits a certain
Killing spinor.  We show that if the hidden symmetry is
tri-holomorphic, then this is equivalent to requiring symmetry along a
higher time and the hidden symmetry determines a `twistor group'
action as introduced by Bielawski \cite{B00}.  This leads to a
construction for the solution to the hierarchy in terms of linear
equations and variants of the generalised Legendre transform for the
hyper-K\"ahler metric itself given by Ivanov \& Rocek \cite{IR96}.  We
show that the ALE spaces are examples of hyper-K\"ahler metrics
admitting three tri-holomorphic Killing spinors.  
These metrics are in this sense
analogous to the 'finite gap' solutions in soliton theory.
Finally we extend the concept of a hierarchy from that of
\cite{DM00} for the four-dimensional hyper-K\"ahler equations to a
generalisation of the conformal anti-self-duality equations and
briefly discuss hidden symmetries for these equations.
}

\newpage
\section{Introduction}
\label{sectFG}
It is well known that finding exact solutions to a nonlinear partial
differential equation (PDE) is greatly simplified by the existence of
symmetries.  In differential geometric language the symmetries of a
hyper-K\"ahler structure or, more generally, anti-self-dual conformal
structures in four dimensions, correspond to (conformal) Killing
vectors. Equations for hyper-K{\"a}hler four-manifolds with conformal
Killing vectors have been studied and, in many cases, solved
\cite{FP79,TW79,DMT00,DT01}.

Apart from natural Lie-point symmetries, integrable soliton
equations possess infinitely many hidden symmetries, which can also be
effectively used to construct solutions.  It is less well know how to
find such solutions in hyper-K{\"a}hler geometry (although we will see
that such solutions have indeed been found in another guise,
\cite{B00, LR1, IR96}). In this paper we shall show that hidden symmetries
correspond to Killing tensors and spinors (which, classically, occur
in Riemannian geometry as additional integrals of the geodesic flow)
and propose two methods of finding hyper-K{\"a}hler metrics with such
symmetries.  

We start by briefly  reviewing a beautiful construction of 
Novikov \cite{N74} which we shall posit as a motivation and
a guiding principle.
Consider the Korteveg--de Vries (KdV) equation 
\be
\label{KdV}
u_{t_1}=6uu_x-u_{xxx},\qquad\mbox{where}\qquad u=u(x, t_1),
\ee
together with the 
associated hierarchy of equations for $u(x, t_1, t_2, ...)$
 \be
\label{kdvhier}
\frac{\p u}{\p t_i}=\frac{\p}{\p x}\frac{\delta H_i}{\delta u}.
\ee
Here
\[
H_0=\int \frac{1}{2}u^2 \d x ,\qquad H_1=\int 
\Big(\frac{1}{2}{u_x}^2+u^3\Big) \d x,  
\qquad H_2=\int\Big(
\frac{1}{2}{u_{xx}}^2-\frac{5}{2}u^2u_{xx}+\frac{5}{2}u^4\Big)\d x,
\qquad ...
\]
are constants of motions which can be found recursively 
by solving  the Riccati equation.
Imposing a constraint 
\be
\label{kdvode}
\frac{\p u}{\p t_k}+c_1\frac{\p u}{\p t _{k-1}}+...+c_k\frac{\p u}{\p
t_0}=c_0
\ee
reduces (\ref{KdV}) to an ODE. This ODE is a completely integrable
Hamiltonian system with $k$ first integrals in involution.
In the simplest non-trivial case the solution is
\[
x=\int{\frac{\d u}{\sqrt{2u^3+c_1u^2+2c_0u+E}}}.
\] 
In the case of the KdV equation, it is then possible to proceed to
obtain explicit formulae for such solutions in terms of theta functions.

For a general integrable system, hidden symmetries are constructed
systematically by studying a hierarchy of commuting flows associated
to the original equations.  A hidden symmetry is then an explicit
point symmetry of the hierarchy which, in particular, include the
higher flows themselves.

In this paper we shall propose an analogous construction of hidden symmetries
for the hyper-K\"ahler equations in four dimensions and its integrable
generalisations, which include quaternionic structures in $4k$ dimensions.  
Recall that a four-dimensional Riemannian manifold 
$({\cM}, g)$ is hyper-K\"ahler if it admits three K\"ahler 
structures $\Sm\sb I$, $\Sm\sb J$ and $\Sm\sb K$ 
compatible with $g$ 
and such that the endomorphisms $I$, $J$, $K$ given by 
$g(IX,Y)=\Sm\sb I(X,Y)$, etc., satisfy $IJ=K=-JI$.
To impose higher symmetries on this system
one needs to:
\begin{itemize}
\item[{\bf (1)}] reformulate  the hyper-K\"ahler condition on a metric
as an integrable PDE (the {\em heavenly equation}),
\item[{\bf (2)}] construct the associated hierarchy,
\item[{\bf (3)}] look for solutions invariant under the hidden symmetries,
and characterise twistor spaces corresponding to these solutions.
\end{itemize}
Steps {\bf (1)} and {\bf (2)} were taken in \cite{Pl75}, and
\cite{BP85,Ta89,St95,DM00}, respectively.  We shall review the
approach taken in \cite{DM00} in \S\ref{hkhiersec}
which focusses on hierarchies associated to 4-dimensional
hyper-K\"ahler spaces.  This is generalised in \S\ref{gcasdhier} to
give a hierarchy associated to general conformally ASD spaces both in
four dimensions and their higher dimensional generalisations such as
quaternion K\"ahler spaces.  The other sections deal with {\bf (3)}.

In \S\ref{hidden_symmetry} and \S\ref{patching_functions} we discuss
symmetry and hidden symmetry reduction in the context of the
4-dimensional hyper-K\"ahler equations.  In \S\ref{hidden_symmetry} we
first discuss and classify symmetries of solutions to the
hyper-K\"ahler hierarchy.  We use the well known twistor description
\cite{TW79} of the Gibbons-Hawking solution as a guiding example in
our analysis of the case of the hierarchy; we show that when the
symmetry is triholomorphic, the solutions of the reduced equations are
linear and we briefly discuss the corresponding twistor theory.  The
hyper-K\"ahler hierarchy is in particular foliated by four-dimensional
Hyper-K{\"a}hler manifolds that admit a hidden symmetry in a variant
of the Ivanov and Rocek construction \cite{IR96}.

In \S\ref{patching_functions} we discuss hyper-K{\"a}hler spaces with
a hidden symmetry, defined to be a space that embeds into a hierarchy
that has a symmetry.  In the general case, we show that a hidden
symmetry corresponds to the existence of a spinor
$K^{B}_{B_0'...B_k'}$ satisfying \be
\label{Killspinor}
{\nabla^{(A}}_{(A'}K^{B)}_{B_0'...B_k')}=0.
\ee
When the hidden symmetry is `triholomorphic' in an appropriate sense,
we find that $K^{B}_{B_0'...B_k'}=\nabla^B_{(B'_0}L_{B'_1\cdots B'_k)}$
for some spinor field $L_{B'_1\cdots B'_k}$ satisfying  
\be
\label{Kspinor}
{\nabla_{A(A'}L_{B_1'...B_k')}=0.}  \ee A non-constant solution
$L_{A_1'...A_k'}$ to the killing-spinor equation (\ref{Kspinor}) is
said to be a Killing spinor of type $(0, k)$ and is also sometimes
known as a solution to the valence--$(0,k)$ twistor equation.  If
$k=2$, then $K_{AA'}={\nabla_{A}}^{B'}L_{A'B'}$ is a tri-holomorphic
Killing vector of the given hyper-K\"ahler space, and the
corresponding metric is of the Gibbons-Hawking form \cite{GH78}.  

If the metric admits a hyper-K\"ahler hidden symmetry and hence
Killing spinor, then the corresponding twistor space admits a globally
defined twistor function $Q$ homogeneous of degree $k$.  This is
because the existence of a Killing spinor implies that the twistor
space of the hyper-K{\"a}hler hierarchy admits an action of a
Hamiltonian vector field with Hamiltonian $Q$ and factor space
$\OO(k)$. The hyper-K\"ahler twistor space therefore arises as an
affine bundle over $\OO(k)$, described by a cohomology class $f\in
H^1(\CP^1, \OO(2-k))$.  The corresponding space-time can be determined
directly and the construction followed through to give explicit
formulae for a basis of the  self-dual two--forms (and therefore for the
metric).

In Section \ref{ale}, we demonstrate that the asymptotically locally
Euclidean (ALE) spaces constructed by Hitchin \cite{H79} and
Kronheimer \cite{K1,K2} admit three triholomorphic hidden symmetries.
We show that the corresponding twistor spaces are elliptic fibrations
over $\OO(k)$ for some $k$, and the transition functions defining these
bundles can be found in terms of elliptic integrals.

Finally in \S\ref{gcasdhier} we extend the concept of a hierarchy from
that of \cite{DM00} for the four-dimensional hyper-K\"ahler equations
to a generalisation of the conformal anti-self-duality equations and
give a brief discussion of hidden symmetries in this context.  

The two-component spinor notation used in the paper is summarised in
the appendix.


\section{The hyper-K\"ahler hierarchy.}\label{hkhiersec}
Let ${\cM}$ be a complex $4$-manifold equipped with a holomorphic
metric $g$ and compatible volume form $\nu$; we shall refer to this
triple as a {\em space-time}.  For a four-manifold with metric, we
have that $T\cM = \bbS\otimes\tbbS$ where $\bbS$ and $\tbbS$ are the
bundles of self-dual and anti-self-dual spinors respectively each
being rank two complex vector bundles on $\cM$.  The metric connection
necessarily preserves this factorisation, and the hyper-K\"ahler
condition is equivalent to the condition that the induced connection
on $\tbbS$ be flat.  (The $I$, $J$ and $K$ act trivially on the
$\tbbS$ factor and on the $\bbS$ factor by Pauli $\sigma$ matrices.)
This implies that, on translation of the 1-form indices to indices $A,
B, \ldots$ denoting membership of $\bbS$ and $A', B', \ldots$ denoting
membership of $\tbbS$, the curvature has the form
$$
R_{AA'BB'CC'DD'}= \varepsilon_{A'B'}C_{ABCD}\varepsilon_{C'D'}\, ,
\qquad C_{ABCD}=C_{(ABCD)}\, ,
$$
where $\varepsilon_{AB}$ and $\varepsilon_{A'B'}$ are skew, and
$\varepsilon_{AB}\varepsilon_{A'B'}$ induces the metric under $T\cM =
\bbS\otimes\tbbS$.  In four-dimensions this amounts to the ASD vacuum
equations \be
\label{ASDvacuum}
\Phi_{ABA'B'}=0\qquad R=0\qquad C_{A'B'C'D'}=0, 
\ee 
Here $R$ is the Ricci scalar, $\Phi_{ABA'B'}$ is the trace-free part
of the Ricci tensor, and $C_{A'B'C'D'}$ is the SD part of the Weyl
tensor (we use the conventions of Penrose and Rindler \cite{PR86}).



We now show how the geometrical characterisation of the hyper-K\"ahler
equations and its hierarchy can be reduced to a differential
equations.
We use a potential
formulation, due to Pleba\'nski \cite{Pl75}, based on the fact that
the equations locally imply the existence of a complex-valued function
$\Th$ and coordinate system $(w,z,x,y)$ such that the metric is given
by \be
\label{metric}
g=2\d w\d x+2\d z\d y-2\Th_{xx}\d z^2-2\Th_{yy}\d w^2+4\Th_{xy}\d w\d z,
\ee
and $\Th$ satisfies so called second heavenly equation
\be
\label{secondeq}
\Th _{xw} +\Th _{yz} + \Th _{xx}\Th _{yy}-{\Th _{xy}}^2=0.
\ee
The associated hierarchy is a differential equation with higher times
generalising one of the formulations of the hyper-K\"ahler equations
in terms of potentials.
We introduce the coordinates $x^{Ai}, A=0,1, i=0...n$ on a $(2n+2)$
dimensional
manifold ${\cal N}$.  The dependent
variable $\Th(x^{Ai})$ satisfies the equations \be
\label{hier2}
\p_{Ai}\p_{Bj-1}\Th -\p_{Bj}\p_{Ai-1}\Th +\{ \p_{Ai-1}\Th, \p_{Bj-1}\Th
\}_{yx}=0,\;\;\;\; i,j=1...n.
\ee
Here $\{..., ...\}_{yx}$ is the Poisson bracket with respect to 
the Poisson structure $\p/\p x^{A0}\wedge\p/\p x_{A0}$. (In order to
make contact with the above for $n=1$, put $x^{A0}=(y,-x)$ and
$x^{A1}=(w,z)$ and note that (\ref{secondeq}) is (\ref{hier2}) with
$i=j=1$.)

This hierarchy has a Lax representation
\be\label{laxsys}
L_{Ai}\Phi= \p_{Ai-1} \Phi + \lambda (\p_{Ai}\Phi+\{\p_{Ai-1}\Theta, \Phi\} )=0
\ee
where $A=0,1$, $i=1,\ldots n$, $\lambda$ is an affine coordinate on
$\CP^1$ and $\Phi(x^{Ai},\lambda)$ is a function on
${\cal N}\times \CP^1$.  It is clear that this provides the
point of contact with the abstract definition.


It is clear from the form of the equations (\ref{laxsys}) that the
space ${\cal N}$ is foliated by 4-dimensional hyper-K\"ahler manifolds
parametrized by $x^{Ai}=$ constant for $i\geq 2$.

In \cite{DM00} the hierarchies were obtained both via a rescursion
operator construction and via a twistor construction.  We shall
summarize these constructions in the remaining part of this section.

Let 
\[
\square_{\Th}=
\p_x\p_w+\p_y\p_z+
\Th_{yy}{\p_x}^2 + \Th_{xx}{\p_y}^2-2\Th_{xy}\p_x \p_y
 \]
denote the wave operator on the ASD
background determined by $\Th$, and let ${\cal W}_{\Th}=$Ker$\;\square_{\Th}$.
\begin{prop}{\em\cite{DM00}}
\label{Rekurencja}
\begin{enumerate}
\item[(i)] Elements of  ${\cal W}_{\Th}$ can be identified with linearised
solutions $\delta\Th$ (i.e. $(\Th+\delta\Th)$ satisfies 
{\em(\ref{secondeq})\em} up to the linear terms in $\delta\Theta${\em)}
of the heavenly equation {\em(\ref{secondeq})\em}.
\item[(ii)] Let $(\delta\Th_1\delta\Th_2) \in{\cal W}_{\Th}\times
  {\cal W}_\Th$. The `recursion operator' ${\cal R}$ is defined to be
  the subspace  ${\cal R}\subset{\cal W}_{\Th}\times {\cal W}_{\Th}$
  on which 
\be
\label{hierdefr}
\p_y (\delta\Theta_2)=(\p_w-\Th_{xy}\p_y+\Th_{yy} \p_x)\delta\Theta_1
,\;\;\; -\p_x(\delta\Theta_2)=(\p_z+\Th_{xx}\p_y-\Th_{xy}
\p_x)\delta\Theta_1 \ee  
\end{enumerate}
\end{prop}
Note that the recursion operator is only an operator in the usual
sense when the subspace
${\cal R}$ can be realized as a graph of a genuine operator 
$R:{\cal W}_{\Th}\longrightarrow {\cal W}_{\Th}$
given by the recursion relations 
\be
\label{hierdef}
R\delta\Theta=\p_y^{-1} ((\p_w-\Th_{xy}\p_y+\Th_{yy} \p_x)\delta\Theta)
,\;\;\; R\delta\Theta=-\p_x^{-1}((\p_z+\Th_{xx}\p_y-\Th_{xy}
\p_x)\delta\Theta ).  
\ee 
This identification with a genuine operator
will only be possible when we impose appropriate boundary conditions.
However, for the definition of the hierarchy as a local system of
equations, we will need only the definition of ${\cal R}$ above.

The first few iterations can be explicitly integrated to give
\[
w\longrightarrow
y\longrightarrow -\Th_x\longrightarrow
\Th_z\longrightarrow...,\qquad
z\longrightarrow-x\longrightarrow-\Th_y\longrightarrow-\Th_w 
\longrightarrow...\;.
\]
We introduce the new coordinates 
$x^{Ai}, A=0,1, i=0...n$. For $i=0, 1, x^{Ai}=x^{AA'}=(w, z, x, y)$ 
are coordinates on ${\cal M}$,
and for $1<i \leq n, x^{Ai}$ are the parameters for
the new flows (with $2n-2$ dimensional parameter space $\X$).
The propagation of $\Th$ along these
parameters is determined by the recursion relations 
$\p_{Ai+1}\Th=R\p_{Ai}\Th$.
However the consistency conditions
imply that  in addition $\Th$ satisfies  the equations (\ref{hier2}).
with Lax system generated by the operators
(\ref{laxsys}).

\smallskip

The twistor theory is summarised in the following:

\begin{theo}{\em\cite{DM00}}
\label{suptw}
There is a 1-1 correspondence between solutions to (\ref{hier2}) on
${\cal M}\times \X$ and
twistor spaces ${\cal PT}_n$ defined as follows.

The twistor space ${\cal PT}_n$ is a three-dimensional complex
manifold with the following structures
\begin{itemize}
\item[1)] a projection $\mu :{\cal PT}\longrightarrow \CP^1$,
\item[2)] a section $s:\CP^1\mapsto {\cal PT}$ of $\mu$
  with normal bundle ${\cal O}(n)\oplus{\cal O}(n)$,
\item[3)] a non-degenerate 2-form $\Sm$
  on the fibres of $\mu$, with values in the pullback from $\CP^1$ of
  ${\cal O}(2n)$.
\item[4)] The choice of coordinate systems and potential $\Th$ in the
second Plebanski form for the hierarchy corresponds on the twistor
space to a choice of point $[o^{A'}]\in\CP^1$ and canonical homogeneous
degree $n$ coordinates $\omega^A$ (i.e., $\Sigma=\d
\omega^0\wedge\d\omega^1$) on a neighbourhood of the fibre of ${\cal
PT}_n$ over $[o^{A'}]$ defined up to $2n^{\mathrm{th}}$ order away from this
fibre.
\end{itemize}
\end{theo}

\noindent
Briefly, the space ${\cal M}\times \X$ is reconstructed as 
the moduli space ${\cal N}$  of 
deformations of the section $s$ given in condition
({\em 2}) above. Then ${\cal N}$ is $2n+2 $ dimensional and we can introduce
coordinates and the function $\Theta$ as follows:

We use homogeneous coordinates $\pi_{A'}=(\pi_{0'},\pi_{1'})$ and
affine coordinate $\lambda=\pi_{0'}/\pi_{1'}$ on $\CP^1$ so that the
point $o$ is represented by $o_{A'}=(0,1)$, or $\lambda=0$.  The
homogeneous coordinates $\omega^A$ (i.e., $(\omega^A,\pi_{A'})\simeq
(c^n\omega^A,c\pi_{A'})$ for $c\in \{\C-0\}$) can be pulled back to ${\cal
N}\times \CP^1$ and the expansion $\omega^A=\pi_{1'}^n\sum_{i=0}^nx^{A,i}\ll^i+
O(\lambda^{n+1})$ and this defines coordinates $x^{A,i}$ on
${\cal N}$.  

Expanding $\omega^A$ further, we discover, as shown in  
 \cite{DM00}, that the twistor coordinates $\om^{A}$ 
pulled back to the correspondence space $\CP^1\times{\cal N}$ can be 
expanded further as 
\begin{eqnarray}\label{twonexp}
{\om_n}^0=(\pi_{1'})^n[x^{0n}+\ll x^{0n-1}+...+\ll^n x^{00}+\ll^{n+1}
\frac{\p\Th}{\p x^{10}}+ \ll^{n+2} \frac{\p\Th}{\p
x^{11}}+...+\ll^{2n+1}\frac{\p\Th}{\p x^{1n}}+...]\;,\nonumber \\
{\om_n}^1=(\pi_{1'})^n[x^{1n}+\ll x^{1n-1}+...+\ll^n x^{10}+
\ll^{n+1}\frac{\p\Th}{\p x^{00}}+\ll^{n+2} \frac{\p\Th}{\p
x^{01}}...+\ll^{2n+1}\frac{\p\Th}{\p x^{0n}}+...]  \;,\;\quad
\end{eqnarray} and this determines $\Th(x^{Ai})$ (up to a constant)
satisfying (\ref{hier2}).  The form of the Lax system (\ref{laxsys})is
determined by the fact that $\omega^B$ are solutions to $L_{Ai}\omega^B=0$.

The form of the expansions (\ref{twonexp}) and equations (\ref{hier2})
can be obtained from the fact that the expansion of
$\Sigma=\d\omega^0\wedge\d \omega^1$ on ${\cal N}\times \CP^1$ in
powers of $\ll$ must truncate after $\ll^{2n}$.

There is a $2n$-dimensional distribution on the `spin bundle'
$D\subset T({\cal N}\times \CP^1)$ that is tangent to the fibres of
the projection ${\cal N}\times\CP^1\longrightarrow {\cal PT}_n$. The distribution $D$
has an identification with ${\cal O}(-1)\otimes \C^{2n}$
and is generated by the Lax system (\ref{laxsys}).

This correspondence is stable under small perturbations of the complex
structure on ${\cal PT}_n$ preserving (1-3).  

One can find the twistor spaces for the four-dimensional
hyper-K\"ahler slices given by $x^{Ai}= const., i\geq 2$ by taking a
sequence of $n-1$ blowups of points in the fibre over $o_{A'}\in
\CP^1$, the choice of point in the fibre to blow up at the $n-i+1$'th
blowup corresponding precisely to the choice of the values of
$x^{Ai}$.  If one wishes to respect Euclidean reality conditions, one
can instead blow up complex conjugate points in the fibres over
$o_{A'}$ and $\hat{o}_{A'}$.


\section{Symmetries of Hyper-K\"ahler hierarchies}\label{hidden_symmetry}
\noindent
For a hyper-K\"ahler space, we can characterise conformal symmetries as
follows: let  $\Sm^{A'B'}=(\Sm^{0'0'}, \Sm^{0'1'}, \Sm^{1'1'})$ be a basis of 
SD two-forms, and let $\Sm=\Sm^{A'B'}\pi_{A'}\pi_{B'}$.
\begin{defi}
A solution to the hyper-K\"ahler equations admits a symmetry if there
exists a vector field $K$ on ${\cal N}$ together with a lift $\tilde{K}$ to
$S^{A'}$ over ${\cal N}$ such that $\Lie_{\tilde{K}}\Sigma=0$.  
\end{defi}
Here the lift must be $\tilde{K}=K+\phi^{A'}_{B'}\pi_{A'}\p/\p\pi_{B'}$
where $\phi^{A'}_{B'} =
-\half\nabla_{BB'}K^{BA'}-\half\nabla_cK^c\varepsilon^{A'}_{B'}$
according to the standard theory of Lie derivatives of spinors, see
\cite{PR86}.     

We will use this also as a definition for symmetries of the
hyper-K\"ahler hierarchy where now $\Sigma$ will be the pull-back from
twistor space to the spin bundle of the corresponding 2-form that is
homogeneous of degree $2n$ in $\pi_{A'}$.
\begin{defi}
A solution to the hyper-K\"ahler hierarchy admits a symmetry if there
exists a vector field $K$ on ${\cal N}$ together with a lift $\tilde{K}$ to
$S^{A'}$ over ${\cal N}$ such that $\Lie_{\tilde{K}}\Sigma=0$.  
\end{defi}
Again, by homogeneity, we must have that the lift must be
$\tilde{K}=K+\phi^{A'}_{B'}\pi_{A'}\p/\p\pi_{B'}$ where
$\phi^{A'}_{B'}$ will be determined by $K$.

In particular, $\tilde{K}$ is therefore in involution with the Lax
distribution $D$ so that it projects down to a global holomorphic
vector field ${\cal K}$ on the twistor space ${\cal PT}$.

 
We can classify symmetries according to the extent to which they
preserve the various structures on $S^{A'}$ or on the twistor space.
This is most easily seen by examining the vertical part $\phi=
\phi_{A'}^{B'}\pi_{B'}\p/\p\pi_{A'}=\tilde{K}-K$ (where here $\tilde K$
denotes the horizontal lift of $K$).  The matrix $\phi_{A'}^{B'}$
generically has two constant eigenvalues (the space-time must be of Petrov 
type $N$ for non-constant eigenvalues to be admissable).
\begin{enumerate}
\item
$K$ will be said to be {\em tri-holomorphic} if the eigenvalues of
$\phi_{A'}^{B'}$ are equal.  The projection of
$\tilde{K}$ to the projective spin bundle, $P S^{A'}$, is then horizontal,
and the vertical part $\phi$ of $\tilde{K}$ on $S^{A'}$ is a multiple
of the Euler homogeneity operator $\pi^{A'}\p/\p \pi^{A'}$, i.e.,
$\phi_{A'}^{B'} =\mu \delta_{A'}^{B'}$.  
Equivalently, it is tri-holomorphic if ${\cal K}$ is tangent to the
fibres of the projection $p:{\cal PT}_n\mapsto\CP^1$.
\item
$K$ is said to be Killing if the trace of $\phi$ vanishes,
$\phi_{A'}^{A'}=0$.  Then $\phi$ 
preserves the form $\pi_{A'}d\pi^{A'}$ on the spin bundle.
\item $K$ is said to be a homothety if the trace of $\phi$,
$\phi_{A'}^{A'}$ is constant, i.e.,
${\Lie}_{\tilde{K}}\varepsilon^{A'B'}\pi_{A'}\d \pi_{B'} =\mu
\varepsilon^{A'B'}\pi_{A'}\d \pi_{B'}$ for some constant $\mu\neq 0$.
 
\end{enumerate}
In the case that the symmetry is not tri-holomorphic, we can further
distinguish the case where $\phi$, on projection to the projective
spin bundle, has one or two zeroes.  In the single zero case in
particular, $\phi_{A'}^{B'}$ will not be diagonalizable.  We will not
pursue this distinction here but see \cite{DMT00,D02}  for a study of such
symmetries on $(++--)$ hyper-K\"ahler space.


\subsection{Examples with triholomorphic Killing symmetry}
We first consider well known reductions of the hyper-K\"ahler
equations, and then analogous reductions of the hierarchy.

\subsubsection{Gibbons-Hawking metrics revisited.}
\label{GHrevisited}
The heavenly equation 
(\ref{secondeq}) with $R(\Th_x)=\Th_z=0$
can be expressed as
\be
\label{ideal}
\d \Th_x\wedge\d x\wedge\d y+\d w\wedge\d \Th_x\wedge\d \Th_y=0.
\ee
Introduce $p:=\Th_x$ and perform a Legendre transform
\[
F(p, y, w):=px(w, y, p) )-\Th(w, z, y, x(w, y, p)).
\]
Then $x=F_p$, $\Theta_y=-F_y$ and (\ref{ideal}) yields the wave
equation \cite{FP79} 
\be
\label{wave}
F_{pw}+F_{yy}=0.
\ee
Implicit differentiation gives
\[
\Th_{yy}=-F_{yy}+\frac{F_{py}}{F_{pp}},\qquad 
\Th_{xy}=-\frac{F_{py}}{F_{pp}},\qquad \Th_{xx}=\frac{1}{F_{pp}},
\]
and so (with the help of (\ref{metric}) and (\ref{wave}))
\begin{eqnarray}
\label{GHmetric}
g&=&F_{pp}(\frac{1}{4}\d y^2+\d w\d p)-\frac{1}{F_{pp}}
(\d z-\frac{F_{pp}}{2}\d y+F_{py}\d w)^2\nonumber\\
&=&\psi(\frac{1}{4}\d y^2+\d w\d p)-\psi^{-1}(\d z+\Om)^2,
\end{eqnarray}
where  $\psi=F_{pp}$ and $\Om=F_{py}\d w-(F_{pp}/2)\d y$ satisfy the
monopole equation $\ast\d \psi=\d \Om$ from (\ref{wave}).
Thus (\ref{GHmetric})
is of the Gibbons-Hawking form \cite{GH78}.

The twistor description is as follows: the twistor
coordinates pull back to the spin bundle as
\begin{eqnarray}
\label{tfunctions}
{\om}^0&=&\pi_{1'}[w+\lambda
y-\lambda^2\Th_x+ \lambda^3\Th_z +...]\;\;,\nonumber\\
{\om}^1&=&\pi_{1'}[z-\lambda x-\lambda^2\Th_y-\lambda^3\Th_w +...]\;.
\end{eqnarray}
The vanishing of $\Theta_z$ implies that the whole series for
$\omega^0$ truncates at 2nd order.  Thus the twistor space admits a
global section of $\OO(2)$, and this is the Hamiltonian with respect to
$\Sigma$, for the holomorphic vector field corresponding to the
Killing field $\p_z=K^{AA'}\p_{AA'}$.  Conversely, given a
tri-holomorphic symmetry, the tri-holomorphicity condition means that
its lift to the spin bundle ${\cal M}$ is horizontal and so on twistor
space, the corresponding holomorphic vector field is tangent to the
fibres of $\mu$.  It also preserves $\Sigma$ and so is Hamiltonian
with Hamiltonian given by a homogeneity degree-2 global function.  We
can choose $\omega^0$ to be this preferred section divided by
$\pi_{1'}$ so that the series for $\omega^0$ terminates after
$\lambda^2$.

Substituting the Legendre transform into (\ref{tfunctions}) yields
\begin{eqnarray}
{\om}^0&=&\pi_{1'}[w+\lambda y-\lambda^2p],\\
{\om}^1&=&\pi_{1'}[z-\lambda F_p+\lambda^2F_y+\lambda^3F_w +...]\;.
\end{eqnarray}
With the definition $\Sm=\d \om^0\wedge\d \om^1|_{\ll=const}$ 
the equation
(\ref{wave}) can be rewritten as $\Sm\wedge\Sm=0.$
The basis of SD two forms can be read off from 
$\Sm=\Sm^{A'B'}\pi_{A'}\pi_{B'}$:
\[
\Sm^{0'0'}=-\d z\wedge\d p+\d y\wedge \d F_p-\d w\wedge\d F_y,\qquad
\Sm^{0'1'}=\d z\wedge\d y +\d w\wedge \d F_p,\qquad
\Sm^{1'1'}=\d z\wedge\d w.
\]
and these determine the metric above.

\subsubsection{Triholomorphic symmetry reductions of the hierarchy}
In this sub-section we shall generalize the construction of
Gibbons-Hawking metrics described in the last subsection, and generate
solutions to the hyper-K{\"a}hler hierarchy such that
$R^n\Th_x=\p_{1n}\Th=0$.  These are the cases of a tri-holomorphic
Killing symmetry.
\begin{prop}
\label{method1} The
Hyper-K\"ahler hierarchy (\ref{hier2}) with symmetry 
$\p\Th/\p x^{1n}=0$ reduces (in appropriate coordinates) to an 
overdetermined system of $n(2n-1)$ linear equations for $F(t^0, ...,
t^{2n})$:
\be
\label{wavehierarchy}
\frac{\p^2F}{\p t^{i+1}\p t^j}=\frac{\p^2F}{\p t^{i}\p t^{j+1}},\qquad
i,j=0,...2n-1.
\ee
\end{prop}
{\bf Proof.} 
Let ${\cal PT}_n$ be the twistor space from Proposition \ref{suptw},
and let  $({\om_n}^A, \pi_{A'})$ be homogeneous 
coordinates\footnote{In the previous section 
${\cal PT}:={\cal PT}_1$ and  $\om^A:={\om_1}^A$
correspond to the standard situation of the nonlinear graviton
construction.} on the neighborhood of $\pi_{A'}=o_{A'}$.

Now impose the symmetry condition, i.e., assume that
$R^n\Th_x=\p_{1n}\Th=0$. Again, the vanising of $\p_{1n}\Theta$
implies that the series (\ref{twonexp}) in $\ll$ for $\omega^0$
truncates at degree $2n$.  [Thus, $\pi_{1'}^n\omega^0$ is a global
holomorphic function of homogeneity degree-$2n$ on ${\cal PT}_n$.
Conversely, again, this symmetry corresponds to a global holomorphic
vector field on ${\cal PT}_n$ that is vertical up the fibres of $\mu$
and preserves $\Sigma$.  It therefore is generated by a global
Hamiltonian $Q$ homogeneous of degree $2n$, and we can take as before
$\omega^0=Q/\pi_{1'}^{n}$.]

We can now perform the Legendre transform
\be
\label{LLegengre}
p^i=\p_{1i}\Th\;\;\; i=0,...,n-1,\qquad 
F(p^j, x^{0j})=\sum_{i=0}^{n-1}p^ix^{1i}(p^j, x^{0j})-
\Th(x^{0j}, x^{1i}(p^j, x^{0j})),\qquad x^{1n}=T.
\ee
Therefore $\p_{0i}F=-\p_{0i}\Th, \p_{p_i}F=x^{1i}$. Define 2n+1 
functions $(t^0, ..., t^{2n})$ by 
\[
t^{n-i-1}=p^i, i=0,...,n-1 \qquad t^{n+i}=x^{0i}, i=0,...n.
\]
This implies
\begin{eqnarray} 
\label{expansion}
{\om_n}^0&=&(\pi_{1'})^n[t^{2n}+\ll t^{2n-1}+...+\ll^{2n}t^0]\nonumber\\
{\om_n}^1&=&(\pi_{1'})^n[T+\ll\frac{\p F}{\p t^0}+\ll^2\frac{\p F}{\p t^1}+...
+\ll^{2n+1}\frac{\p F}{\p t^{2n}}+...]\;.\\
\end{eqnarray}
The equations 
(\ref{wavehierarchy}) arise from the vanishing of coefficient
$\ll^{2n+2}$ in 
$\d {\om_n}^0\wedge\d {\om_n}^1$:
\[
\sum_{i=0}^{2n-1}\d t^i\wedge\d\frac{\p F}{\p t^{i+1}}=0.
\]
It can also be verified by cross-differentiating 
that all integrability conditions
for the system (\ref{wavehierarchy})
are satisfied.
\koniec

The geometry on twistor space can be understood as follows. The
section $Q$ of $\OO(2n)$ generates a Hamiltonian flow
\[
{\cal K}=\varepsilon^{AB}\frac{\p Q}{\p \omega^A}\frac{\p}{\p\om^B}\, ,
\quad \varepsilon^{AB}=-\varepsilon^{BA}\, , \;\varepsilon^{01}=1
\]
on the extended twistor space ${\cal PT}_{n}$. This flow corresponds to
${\cal K}^{AA_1'...A_n'}\p_{AA_1'...A_n'}=\p/\p x^{1n}$ on ${\cal N}$.
Since $Q$ is constant along ${\cal K}$, the quotient space 
${\cal PT}_{n}/{\cal K}$ is the total space of
$\OO(2n)\rightarrow\CP^1$ where the map to $\OO(2n)$ is furnished by
$(\omega^A,\pi_{A'})\longrightarrow (Q,\pi_{A'})$. 

The full twistor space ${\cal PT}_n$ is an affine line bundle with
trivial underlying translation bundle and so it corresponds to an
element $G(Q,\pi_{A'})$ of $H^1(\OO(2n),\OO)$.

Sections of $\OO(2n)\rightarrow \CP^1$ are parametrised by $\C^{2n+1}$
with coordinates ${\bf t}=(t^0, ...,t^{2n})$. The $2n+2$ dimensional space of
sections of ${\cal PT}_n\rightarrow \CP^1$ maps onto this with fibre
$\C$ parametrized by $T$.  Choosing linear coordinates up the line
bundle $\eta$ and $\tilde{\eta}$ over open sets $U$ and $\tilde{U}$,
the problem of lifting a curve $L_{\bf t}$ in $\OO(2n)$ to one in ${\cal
PT}_n$ is one of finding a trivialization of the line bundle over
$L_t$, i.e. of finding functions $g({\bf t},\pi_{A'})$ and
$\tilde{g}({\bf t},\pi_{A'})$  such that, on $U\cap \tilde{U}$,
$g-\tilde{g}=G$ on restriction to $L_t$ where $G$ is the log of the
patching function for the line bundle, and therefore has homogeneity
degree zero.  We can then take $\eta=T+g$.
We can give a formula for $g$ as
$$
g({\bf t},\ll)= \oint_\gamma G(t^{2n}+\zeta
t^{2n-1}+...+\zeta^{2n}t^0,\zeta)\d\zeta /(\zeta -\ll) \, .
$$
where the contour $\gamma$ is taken in $L_t\cap U\cap \tilde{U}$
surrounding $\zeta$ in $U$.  With an identical expression for
$\tilde{g}$ but with contour $\tilde{\gamma}$ such that
$\gamma-\tilde{\gamma}$ surrounds $\zeta$, we see that $g-\tilde{g}=G$
follows from Cauchy's integral formula.  Then the expression for 
the expansion of the coordinate $\eta=T+g$ about $\ll=0$ is
$$
\eta=T + \sum_{i=0}\ll^i \oint G(Q(t,\zeta),\zeta)\d \zeta/\zeta^{i+1}\, .
$$
It can then be seen that if we define
$$
F({\bf t})
=\oint_{\Gamma}  \frac{G(t^{2n}+\ll t^{2n-1}+...+\ll^{2n}t^0,\ll)}{\ll^2}
\d\ll \, .
$$
then $F$ clearly satisfies the
equations of $\ref{wavehierarchy}$ and we obtain the  expansions
(\ref{expansion}) for $(\omega^0, \omega^1)=(Q/\pi_{1'}^n,
\pi_{1'}^n \eta)$. These can then be used to obtain concrete
expressions for $\Sigma=\d \omega^0\wedge \d \omega^1$ to determine
the geometric structures of the hyper-K\"Ahler hierarchy. 

Clearly we have:
\begin{lemma}\label{hkhierred}
The full space of the hierarchy is foliated by hyper-K{\"a}hler
four-manifolds with $x^{Ai}=$ constant for $i>1$ and metric 
\be
\label{metric_on_leaf}
2\d x^{10}\d x^{01}+2\d x^{11}\d x^{00}-
2\frac{\p^2 \Th}{\p {(x^{10})}^2}(\d x^{11})^2
-2\frac{\p^2 \Th}{\p {(x^{00})}^2}(\d x^{01})^2
-4\frac{\p^2 \Th}{\p x^{00}\p x^{10} }\d x^{01}\d x^{11}\, .
\ee
\end{lemma}
This gives a variant of the 
Legendre transform of Ivanov and Rocek \cite{IR96}
(see also \cite{B00}).  
\subsection{Example for $n=2$}
We saw above that for $n=1$ the construction is equivalent to the
Gibbons-Hawking anzatz. The $n=2$ case goes as follows:

Let $ F_{i}:={\p_{t^i} F}$. Implicit differentiation of
$\p_{0i}F=-\p_{0i}\Th, \p_{p_i}F=x^{1i}$ with respect to $p_0, p_1, y$
yields
\[
\Th_{xx}=-\frac{F_{00}}{M},\;
\Th_{xy}=\frac{-F_{01}F_{02}+F_{00}F_{12}}{M},\;
\Th_{yy}=
-F_{22}-\frac{F_{00}(F_{12})^2-2F_{01}F_{12}F_{02}+F_{11}(F_{02})^2}{M},
\]
where $M:=(F_{01})^2-F_{00}F_{11}$.  The metric (\ref{metric}) with
$x=x(t^i), z=z(t^i)$ is defined on the surface $F_4=0$. 
The formula for the metric in terms of $F$ is not very illuminating,
but we shall give it for the sake of completeness:
\begin{eqnarray}
g&=&M^{-1}(F_{10}N\d t^1\d t^2+F_{00}N\d t^0\d t^2+F_{11}(F_{01})^2(\d
t^2)^2+F_{11}^3(\d t^3)^2+2F_{01}(F_{11})^2\d t^2\d t^3\nonumber\\
&&+2F_{01}(F_{00})^2\d t^0\d t^1 +F_{00}^3(\d t^0)^2
+F_{00}(F_{01})^2(\d t^1)^2+[F_{11}N+F_{01}F_{00}F_{03}]\d t^1\d t^3\nonumber\\
&&+[3F_{01}F_{00}F_{11}-(F_{01})^3 +(F_{00})^2F_{03}]\d t^0\d t^3),
\end{eqnarray}
where $N:=(F_{01})^2+F_{00}F_{11}$.

\section{
Hyper-K\"ahler spaces with hidden symmetries}
\label{patching_functions}


If we have a hyper-K\"ahler space that embeds into a hyper-K\"ahler
hierachy that admits a symmetry we will say that that the original
hyper-K\"ahler space admits a {\em hidden symmetry}.

The first question we wish to address is of how to recognise
when a hyper-K\"ahler space admits such a hidden symmetry.  

\begin{prop}\label{hiddensym}
If a hyper-K\"ahler space ${\cal M}$ admits a hidden symmetry then it
admits a solution to the equation 
\be\label{cke}
\nabla_{(A'_1}^{(A}K^{B)}_{A'_2\cdots A'_n)} =0 \ee
\end{prop}

\noindent
{\bf Proof:} This result is most easily seen from the twistor theory.
The symmetry vector $K$ gives rise to a global holomorphic vector
field $\cK$ on the twistor space for $\cN$, $\cPT_n$.  The twistor
space $\cPT_1$ for $\cM$ is a region in the blowup of $\cPT_n$ at a
number of points.  Thus we have a map $p:\cPT_1 \rightarrow \cPT_n$.  So rather than consider $\cK$ itself, we consider
the 2-form $\cK\hook \nu_n$ 
weight $2n+2$ where
$\nu_n\in\Gamma(\cPT_n,\Omega^3(2n+2))$ is given by $\Sigma_n\wedge
\pi_{A'}\d\pi^{A'}$.

This 2-form can be pulled back to give $p^*\cK\hook \nu_n$ a global 2-form of weight
$2n+2$ on the twistor space $\cPT_1$ for $\cM$.  This 2-form can then
be pulled back to give  a 2-form on the spin-bundle $\bbS^{A'}$ which
must take the form
$$
\cK\hook\nu_n= K^{AA'_1\cdots A'_{2n-1}}e_A^{A'_{2n}}\pi_{A'_1}
\cdots \pi_{A'_{2n}}\wedge \pi_{B'}\d\pi^{B'} +
\chi_{A'_1\cdots A'_{2n}}\pi^{A'_1}\cdots \pi^{A'_{2n}}\Sigma_1
$$
for some $K^{AA'_1\cdots A'_{2n-1}}$ and $\chi_{A'_1\cdots A'_{2n}}$.

The condition that this 2-form descends to twistor space is the
condition that 
\[
\pi^{A'}\nabla_{AA'}\hook \d(\cK\hook\nu_n)=0.
\]  
This leads to equation (\ref{cke}) and 
$$
\nabla_{(A'_1}^{A}K_{A'_2\cdots A'_{2n})A} =\chi_{A'_1\cdots
A'_{2n}}\, , \qquad \nabla_{A(A'}\chi_{A'_1\cdots
A'_{2n})}=0\, .
$$
However, it can be checked that these two equations are a
consequence of (\ref{cke}), if the first equation is taken to be the
definition of $\chi_{A'_1\cdots A'_{2n}}$.
\koniec

\subsection{The case of a hidden tri-holomorphic Killing symmetry}
The case of a hidden tri-holomorphic Killing symmetry reduces to
linear equations and can be worked through completely modulo some
intergations and solving for implicit functions.  This is effectively
the case studied by Ivanov and Rocek \cite{IR96} and generalised by
Bielawski \cite{B00}.

In the tri-holomorphic Killing case, we have

\begin{lemma}\label{kstf}
Suppose $\cM$ admits a hidden tri-holomorphic Killing symmetry, then
$\chi_{A'_1\cdots A'_{2n}}=0$ and there exists  a spinor
$\phi_{A'_2\cdots A'_{2n}}$ such that 
$$
\nabla_{AA'}\phi_{A'_1\cdots A'_{2n}}=
K_{A(A'_1\cdots A'_{2n-1}}\varepsilon_{A'_{2n})A'}
$$
\end{lemma}

\noindent
{\bf Proof:}
The vanishing of $\chi_{A'_1\cdots A'_{2n}}=0$ follows from the fact
that $\cK$ is tangent to the fibres of twistor space over $\CP^1$.
The existence of $\phi_{A'_1\cdots A'_{2n}}$ follows from the fact
that $\cK$ is Hamiltonian with respect to the symplectic forms
$\Sigma_n$ up the fibres of $\mu$ and so is generated by a Hamiltonian
$Q\in\Gamma(\OO(2n))$.  On pullback to the spin bundle
$Q=\phi_{A'_1\cdots A'_{2n}}\pi^{A'_1}\cdots \pi^{A'_{2n}}$ and the
condition that $Q$ descends to twistor space is $\pi^{A'}\nabla_{AA'}
Q=0$ and this gives the equation above.
\koniec

The above lemma shows that a Killing spinor on an ASD vacuum
determines a function $Q$ homogeneous of degree $k$ on its twistor
space ${\cal PT}$. This in turn implies 
\begin{lemma} If an ASD vacuum space-time admits a Killing spinor, its
  twistor space ${\cal PT}$ is an affine line bundle over $\OO(k)$ with
  underlying translation bundle $\OO(2-k)$.  
\end{lemma}

\noindent
{\bf Proof:} 
The existence of a global twistor function homogeneous
of degree $k$ on ${\cal PT}$ gives a projection onto $p:{\cal
  PT}\rightarrow \OO(k)$. Furthermore the fibre is spanned by the
Hamiltonian vector field of Q with respect to $\Sigma$, in local
coordinates,
$$
{\cal K}=\varepsilon^{AB}\frac{\p Q}{\p\omega^A}\frac{\p}{\p
  \omega^B}\, .
$$
This is a vector field with values in $O(k-2)$.  This gives each fibre
an affine linear structure which is twisted globally by $\OO(2-k)$,
since, if $a$ is a local section of $\OO(2-k)$ over $\CP^1$, then
$a{\cal K}$ is a vector whose flows determine an action of $\C$.  Thus
${\cal PT}\rightarrow\OO(k)$ is an affine line bundle over $\OO(k)$ with
underlying translation bundle $\OO(2-k)$.  \koniec

Such affine line bundles are classified by elements $[f]$ of $H^1(\OO(k),
\OO(2-k))$.  In a Cech description, cover $\OO(k)$ by open sets, $U_i$,
and represent $[f]$ by its Cech representative $f_{ij}\in\Gamma(\OO(2-k),
U_i\cap U_j)$.  Then ${\cal PT}$ is constructed by patching together
the total space of $\OO(2-k)\rightarrow U_i$ to $\OO(2-k)\rightarrow
U_j$ by translating the zero section by $f_{ij}$.  The data $[f]$
therefore determines the twistor space.  This proves the first part of:


\begin{theo} \label{maintheorem}
There is a one-to-one correspondence between ASD vacuum space-times
$({\cal M},g)$ admitting a valence $(0,k)$ Killing spinor and elements
$[f]$ of $H^1(\OO(k), \OO(2-k))$.

\label{method2}
In this case, for $k\geq 3$, $\cM$ admits a natural map into
$\C^{k+1}=\odot^k\bbS^{A'}$ which we coordinatise with $t^{A'_1\cdots
A'_k}$.   The hyper-k\"ahler space
${\cal M}$
 is determined as a subset of $\C^{k+1}$ by the $k-3$ constraints
\[
f_{A'_1\cdots
A'_{k-4}}:=\oint_{\Gamma}\pi_{A_1'}...\pi_{A_{k-4}'}f(Q,
\pi_{A'})\pi\cdot\d\pi=0.  
\]

The basis of $SD$ two forms for $g$ is then given by the restriction
of the forms
\be
\label{SDtwoform}
\Sm^{A_{1}'B_{1}'}=
\psi_{B_2'...B_{k}'A_2'...A_{k-2}'}\d t^{A_1'...A_k'}\wedge
\d {t^{B_2'...B_{k}'}}_{A_k'},
\ee
to $\cM$, where
\be
\label{the_field}
\psi_{A_1'...A_{2k-4}'}=\frac{1}{2\pi i}\oint_{\Gamma}
\pi_{A_{1}'}...\pi_{A_{2k-4}'}
\frac{\p f}{\p Q}\pi\cdot\d\pi
\ee
is a field determined by an arbitrary element of $H^1(\OO(k),
\OO(2-k))$. 

\end{theo}

\noindent
{\bf Proof.}  If one wishes to obtain the space-time $({\cal M}, g,
\nu)$ determined by a given a twistor space, the first task is to
locate the 4-dimensional family of sections of the fibration ${\cal
PT}\rightarrow \CP^1$.

Let $t^{A_1'...A_k'}=t^{(A_1'...A_k')}$ be coordinates on the
$\C^{k+1}=\odot^k {\mathbb S}^{A'}$ parameter space of sections
$\sigma_t: \pi_{A'}\rightarrow
Q=t^{A_1'...A_k'}\pi_{A_1'}...\pi_{A_k'}\in \Gamma(\OO(k))$.  Sections
of ${\cal PT}\rightarrow \CP^1$ determine sections of $\OO(k)$ by
projection onto $\OO(k)$.  However, the affine line bundle ${\cal PT}$
only admits a section over some $\sigma_t$ if the cohomology class
$[f]$ vanishes on restriction to $\sigma_t$.  If $[f]$ vanishes
on restriction to $\sigma_t$, ${\cal PT}$ restricts to become the line
bundle $\OO(2-k)$ so that there is a $3-k$-dimensional family of
sections over $\sigma_t$ for $3-k>0$ or just the $0$-section
otherwise.

To obtain explicit formulae,  we first note that $[f]$ determines a
field 
$$
f_{A'_1\cdots A'_{k-4}}(t^{B'_1\cdots B'_k})
=\oint_{\gamma\subset \sigma_t}
\pi_{A'_1}\cdots 
\pi_{A'_{k-4}} f\pi_{B'}\d \pi^{B'}
$$
on $\C^{k+1}$ (here we express the natural pairing as a contour
integral over some contour $\gamma_t$ in $\sigma_t$ where here $\chi$
is a Cech representative).  This vanishes at some $t^{A'_1\cdots
A'_k}$ iff $[f]$ vanishes on the corresponding $\sigma_t$.  Thus,
for $3-k>0$, ${\cal M}$ is fibred over the zero set of
$\chi_{A'_1\cdots A'_{k-4}}$ in $\C^{k+1}$ with $3-k$-dimensional
fibres, and is simply identified with this zero set for $k\geq 3$.

In order to calculate the SD 2-forms associated to the space-time, we
use the method of Gindikin \cite{G82}, and pullback $\Sigma$ to the
spin-bundle. To this end, introduce local homogeneous coordinates
$(\pi_{A'}, Q, \zeta_i)= (c \pi_{A'}, c Q, c^{2-k}\zeta_i)$
on each set $U_i$ of some Stein cover twistor space; here $\zeta_i$ is a
fibre coordinate up the fibres of the affine line bundle ${\cal
PT}\rightarrow \OO(k)$ on $U_i$ with patching relations
$\zeta_i=\zeta_j+f_{ij}$ on $U_i\cap U_j$.  In these coordinates
$$
\Sigma=\d_h Q\wedge\d_h\zeta_i\, ,$$ where $\d_h$ denotes the exterior
derivative in which $\pi_{A'}$ is held constant, i.e., horizontal on the
spin bundle over space-time (although slightly confusingly, vertical
with respect to the fibration ${\cal PT}\rightarrow \CP^1$).  This form
$\Sigma$ is globally defined on vector fields tangent to the fibres of
${\cal PT}\rightarrow \CP^1$ as $f_{ij}$ does not depend on
$\zeta_i$.

In order to evaluate this, we need to find the values of $\zeta_i$ on
the sections of ${\cal PT}\rightarrow\CP^1$.  These are obtained by a
splitting formula due to Sparling.  

On a $\sigma_t$ for which $f_{A'_1\cdots A'_{k-4}}(t^{B'_1\cdots
B'_k})=0$, we can find $\zeta_i(t,\pi_{A'})$ such that
\begin{equation}\label{split}
\zeta_i(t,\pi_{A'})=\zeta_j(t,\pi_{A'})+
f_{ij}(t^{A'_1\cdots A'_k}\pi_{A'_1}\cdots\pi_{A'_k}, \pi_{B'})\, .
\end{equation}
For $k\geq 3$  this solution will be unique, but for $k<3$ we will be
free to add $x^{A'_1\cdots A'_{k-2}}\pi_{A'_1}\cdots \pi_{A'_{k-2}}$
to the solution.

In the formula for $\Sigma$ we can rearrange so that we have
$$
\Sigma=\d_h Q\wedge \d_h\zeta_i =\d
t^{A'_1\cdots A'_k}\wedge\d_h(\pi_{A'_1}\cdots\pi_{A'_k} \zeta_i)\, .
$$
Applying $\d_h$ to equation (\ref{split}), and multiplying by $k-3$ of the
$\pi$s, we obtain
$$
\d_h(\pi_{A'_1}\cdots\pi_{A'_{k-3}}\zeta_i)
=\d_h(\pi_{A'_1}\cdots\pi_{A'_{k-3}} \zeta_j) + \frac{\p f_{ij}}{\p Q}
\pi_{A'_1}\cdots\pi_{A'_{k-3}} \pi_{B'_1}\cdots\pi_{B'_k}
\d t^{B'_1\cdots
B'_k}\, .
$$
The cocycle $\p f_{ij}/\p Q$ defines a class in $H^1(\OO(k),
\OO(2-2k))$ so that the expression above takes values in $\OO(-1)$ on
$\CP^1$ for each fixed $t$.
Thus the splitting as stated exists and is unique since
$H^0(\CP^1,\OO(-1))=H^1(\CP^1,\OO(-1))=0$. 
This gives 
$$
\d_h(\pi_{A'_1}\cdots\pi_{A'_{k-3}}\zeta_i)=k_{i,A'_1\cdots A'_{k-3}
B'_1\cdots B'_k}\d t^{B'_1\cdots
B'_k}\, ,
$$
where $k_{i,A'_1\cdots A'_{2k-3}}$ is defined by the splitting relation
$$
k_{i,A'_1\cdots A'_{2k-3}}-k_{j,A'_1\cdots
A'_{2k-3}}=\pi_{A'_1}\cdots\pi_{A'_{2k-3}}\frac{\p  f_{ij}}{\p Q}
$$
or alternatively the contour integral formula
\be
\label{ijsplit}
k_{i,A_{1'}...{A_{2k-3}'}}=
\oint_{\Gamma_i}
\rho_{A_1'}...\rho_{A_{2k-3}'} 
\frac{1}{\pi\cdot \rho}\frac{\p f_{ij}(Q, \rho_{A'})}{\p
Q}\rho\cdot\d\rho,
\ee
where for simplicity we assume a two set cover and the contour
$\gamma_i$ is chosen so that $\gamma_i - \gamma_j$ surrounds $\pi =
\rho$.   It follows that 
\be\label{psidef}
\psi_{A'_1\cdots A'_{2k-4}}=
\pi^{A'_{2k-3}}k_{A'_1\cdots A'_{2k-4}A'_{2k-3}} = \oint_{\gamma_i}
\rho_{A_1'}...\rho_{A_{2k-4}'}  
\frac{\p f_{ij}(Q, \rho_{A'})}{\p
Q}\rho\cdot\d\rho,
\ee
is the field on $\C^{k+1}$ naturally associated to $\p f/\p Q$.   

We therefore obtain the formula
$$
\Sigma =\d t^{A'_1\cdots A'_k}\wedge  \pi_{A'_1} \pi_{A'_2} 
\pi_{A'_3} k_{A'_4\cdots A'_{2k}}\d t^{A'_{k+1}\cdots A'_{2k}}\, .
$$ 
Define the indexed 2-forms $\Sigma^{(B'_1 B'_2)| (A'_5\cdots A'_{2k})}$
by\footnote{This indexed 2-form can be represented as
$$
\Sigma^{(B'_1 B'_2)| (A'_5\cdots A'_{2k})}=\frac{3k}{2}
  \d t^{(B'_1 B'_2)| (A'_5\cdots
A'_{k+1}}_{C'}\wedge \d t^{A'_{k+2}\cdots A'_{2k})C'} 
+ a_k\varepsilon^{(A'_5|(B'_1} \varepsilon^{B'_2)|A'_6} 
\d t^{A'_7\cdots
A'_{k+1}}_{C'D'E'}\wedge \d t^{A'_{k+2}\cdots A'_{2k})C'D'E'},
$$ where 
$a_k$ is a conbinatorial constant depending on $k$.}
$$
\d t^{B'_1\cdots B'_3 (A'_4\cdots A'_k}\wedge \d t^{A'_{k+1}\cdots A'_{2k})}
= \varepsilon^{(A'_4|(B'_1} \Sigma^{B'_2 B'_3)| A'_5\cdots A'_{2k})C'} 
$$
With this, we find that a $\pi_{A'}$ is contracted onto $k_{i,
A'_1\cdots A'_{2k-3}}$ so that (\ref{psidef}) gives
$$
\Sigma=\pi_{A'_1}\pi_{A'_2}\psi_{A'_3\cdots A'_{2k-2}}\Sigma^{A'_1\cdots
A'_{2k-2}} 
$$
Thus, the result follows.\koniec

{\bf Remarks}
\begin{itemize}
\item
If $k>3$ then there exists a potential for
$\psi_{A_1'...A_{2k-4}'}$;
\[
\psi_{A_1'...A_{2k-4}'}=\p_{A_{k-3}'....A_{2k-4}'}f_{A_1'...A_{k-4}'},
\]
where
\[
f_{A_1'...A_{k-4}'}=
\oint_{\Gamma}\rho_{A_1'}...\rho_{A_{k-4}'}f\rho\cdot\d\rho.
\]
The space-time is a four-dimensional  
surface $\chi_{A_1'...A_{k-4}'}=0$ in $k+1$
dimensional moduli space of $\OO(k)$ sections coordinatized by 
$x^{A_1'...A_k'}$.
\item  $k=3$ 
the field $\psi_{A_1'A_2'}$ doesn't have a potential, and 
no conditions have to be imposed to on the moduli space of  $\OO(3)$
to find the space time. This is because 
$H^0(\CP^1, \OO(3))=H^0(\CP^1, \OO(1)\oplus\OO(1))$, and 
$x^{A_1'A_2'A_3'}$ 
has as many components as $x^{AA'}$. 
Ward \cite{W78} regards $\psi_{A_1'A_2'}$ 
as a self-dual Maxwell field on $\C^4$. 
\item
The case $k=2$ implies the existence of a tri-holomorphic Killing
vector and was consider by Tod and Ward in \cite{TW79}. Now
\[
\psi=\oint_{\Gamma}\frac{\p f}{\p Q}\rho\cdot\d \rho
\]
is a solution to the
three-dimensional wave equation.  The relation between our
construction and the description of the Gibbons-Hawking metric 
form Section \ref{GHrevisited} is given by
\[
F(x^{A'B'})=\oint_{\Gamma}\frac{G(\pi_{A'}, Q)}{(\pi\cdot o)^{2}}\pi\cdot\d\pi,\qquad
\psi=F_{pp}=\oint_{\Gamma}(\pi\cdot o)^2\frac{\p^2 G}{\p
  Q^2}\pi\cdot\d\pi,
\]
where
\[ 
x^{A'B'}:=
\left (
\begin{array}{cc}
-p&y/2\\
y/2&w
\end{array}
\right ).
\]
\end{itemize}
\subsection{Relation between the two constructions}
In this section we  relate the hyper-K\"ahler slices (\ref{hkhierred})
of the symmetric hierarchy introduced in
\S\ref{hidden_symmetry} to the method described above.

\begin{prop}
\label{two_methods}
Let $f\in H^1(\OO(2n),\OO(2-2n))$ give rise to the ZRM field
(\ref{the_field})  with $k=2n$. Then hyper-K\"ahler 
metrics arising from Proposition \ref{method1}, 
form a subclass of metrics  from Proposition 
\ref{method2} if
\[
F=\oint_{\Gamma}(\pi\cdot o)^{-2}G(Q,\pi_{A'})\pi\cdot\d \pi,\qquad \mbox{where}\qquad
(\pi\cdot o)^2\frac{\p G}{\p Q}=f\in H^1(\OO(2n),\OO(2-2n)),
\]
where $o_{A'}$ is a constant spinor.
\end{prop}
{\bf Proof.} Let $(Q,\pi_{A'})$ 
be homogeneous coordinates on the total space of
$\OO(2n)$ bundle. Let us choose a constant spinor $o_{A'}$ and 
parameterize a section of the $\OO(2n)\rightarrow
\CP^1$ by $2n+1$ complex numbers
\[
x^{A_1'...A_{2n}'}=\frac{\p^{2n} Q}{\p\pi_{A_1'}...\p\pi_{A_{2n}'}}
|_{\pi_{A_i'}=o_{A'}}.
\]
The coordinates
$x^{A_1'...A_{2n}'}$ on ${\C^{2n+1}}$ 
correspond to $t^{0}, ... , t^{2n}$ by
\[
t^{i}={2n \choose i}x^{A_1'A_2'...A_{2n}'}o_{A_1'}...o_{A_i'}
\iota_{A_{i+1}'}...\iota_{A_{2n}'}(-1)^{n-i},\qquad
i=0,...,2n.
\]
Define
\[ 
\frac {\p}{\p t^i}=\iota^{A_1'}...\iota^{A_{i}'}o^{A_{i+1}'}...o^{A_{2n}'}
\frac{\p}{\p x^{A_1'...A_{2n}'}}.
\]  
Let 
\[
F=\oint_{\Gamma}(\pi\cdot o)^{-2}G(Q,\pi_{A'})\pi\cdot\d \pi,\qquad \mbox{where}\qquad
(\pi\cdot o)^2\frac{\p G}{\p Q}=f\in H^1(\OO(2n),\OO(2-2n)).
\]
We have
\begin{eqnarray*}
f_{A_1'...A_{2n-4}'}&=&\oint_{\Gamma}\pi_{A_1'}...\pi_{A_{2n-4}'}f
\pi\cdot\d\pi=\oint_{\Gamma}\pi_{A_1'}...\pi_{A_{2n-4}'}
(\pi\cdot o)^2\frac{\p G}{\p Q}
\pi\cdot\d\pi\\
&=&o^{A_{2n-3}'}...o^{A_{2n}'}\frac{\p}{\p x^{A_1'...A_{2n}'}}
\oint_{\Gamma}(\pi\cdot o)^{-2}G\pi\cdot\d\pi=\frac{\p F}{\p
  x^{A_1'...A_{2n-4}'0'0'0'0'}} .
\end{eqnarray*}
Therefore fixing $f_{A_1'...A_{2n-4}'}$ is equivalent to
fixing $\p F/\p t^i$ for $i<2n-3$.
Moreover the global twistor function is given by
\[
Q=(\om_n\cdot\iota)(\pi\cdot\iota)^n=(\pi\cdot\iota)^{2n}
\sum_{i=0}^{2n}\ll^it^{2n-i}.
\]
\koniec
\section{ALE spaces revisited; Finite-Gap solutions of ASD vacuum
equations}
\label{ale}
One way to generalise the Novikov construction 
of finite gap solutions of the Korteweg de Vries equation
to hyper-K\"ahler equations would be to 
study solutions to (\ref{secondeq})
which are invariant under 
three commuting hidden symmetries that we shall take to be tri-holomorphic:
\[
\p_{T_1}:=\sum_{i=1}^ka_i \frac{\p }{\p t_{i}},\qquad 
\p_{T_2}:=\sum_{i=1}^lb_i \frac{\p }{\p t_{i}},\qquad
\p_{T_3}:=\sum_{i=1}^mc_i \frac{\p }{\p t_i},\;\;\;\;\;\mbox{where}\;a_i, b_i, c_i \;\mbox{are constant},
\]
and the propagation of $\Th$ along the parameters $t_i$ is determined
by the recursion relations (\ref{hier2}).  This would reduce
(\ref{secondeq}) down to an ODE. 

Rather than performing the explicit reduction to an ODE, we see from
the twistor picture that the twistor space must  have three
projections onto the total space of the line bundle $\OO(n)$ for three
values of $n$.  Thus we have a map of the twistor space $\cPT$ into
$\OO(p)\oplus\OO(q)\oplus\OO(r)$ and so $\cPT$ can be realized as a
hypersurface in this space (although there may need to some blowup or
resolution of singularities where the map fails to be an embedding).
If we realize $\cPT$ as the zero set of a function $F$ taking values
in a line bundle of degree $s$, then we must have, for rational curves
to have the appropriate normal bundle, that $p+q+r=2+s$.

We will now see that the ALE hyper-K\"ahler spaces falls precisely
into this above class.

It is well know that hyper-K\"ahler manifolds $({\cal M}, g)$ which
have the topology of $\R^4$ at infinity, and approach the flat
Euclidean metric $\eta=\d {x_1}^2+...+\d {x_4}^2$ sufficiently fast,
in the sense that \be
\label{ALE}
g_{ab}=\eta_{ab}+O(r^{-4}),\qquad (\p_a)^p(g_{bc})=O(r^{-4-p}),\qquad
r^2=x_1^2+...+x_4^2
\ee
have to be flat. A weaker asymptotic condition one can impose on
$g$ is assymptoticaly locally Euclidean (ALE).

The ALE spaces are non-compact, complete hyper-K\"ahler
manifolds which satisfy (\ref{ALE}) only
locally for $r\rightarrow \infty$. Globally the neighbourhood of
infinity must look like $S^3/\Gamma\times \R$, where $\Gamma$ is a
finite
group of isometries acting freely on $S^3$ (a Kleinian group).
These manifolds belong to the class of {\em gravitational instantons}
because their curvature is localised in a `finite region` of a space-time.

Finite subgroups of $\Gamma\subset SU(2)$ correspond Platonic solids in $\R^3$.
They are the cyclic groups, and the binary dihedral,
tetrahedral, octahedral and icosahedral groups
(one can think about the last three as M\"obius transformations
of $S^2=\CP^1$ which leave the points  corresponding to vertices of a
given Platonic solid fixed). 
Each of them can be
related
to a Dynkin diagram of a simple Lie algebra. 
All Kleinian groups act on $\C^2$, and the `infinity`
$S^3\subset\C^2$. Let $(z_1, z_2)\in\C^2$. 
For each $\Gamma$ there exist three invariants $x, y, z$
which are polynomials in $(z_1, z_2)$ invariant under $\Gamma$.
These invariants satisfy some algebraic relations which 
we list below:
{\small
\begin{center}
\begin{tabular}{p{1cm}|lll}
\multicolumn{4}{c}{}\\
&{\bf Group} & {\bf Dynkin diagram} & {\bf Relation} $F_{\Gamma}(x,y,z)=0$\\
&cyclic            &$A_k$& $xy-z^k=0$\\
&dihedral          &$D_{k-1}$&$ x^2+y^2z+z^{k}=0$\\
&tetrahedral       &$E_6$& $x^2+y^3+z^4=0$\\
&octahedral        &$E_7$& $x^2+y^3+yz^3=0$\\
&icosahedral       &$E_8$& $x^2+y^3+z^5=0$
\end{tabular}
\end{center}
\small}
In each case
\[
\C^2/\Gamma\subset \C^3=\{(x, y, z)\in\C^3, F_{\Gamma}(x, y, z)=0\}
\]
The manifold ${\cal M}$ on which an ALE metric is defined is obtained by
minimally resolving the singularity at the origin of $\C^2/\Gamma$. This
desingularisation is achieved by taking ${\cal M}$ to be the zero
set of
\[
\widetilde{F}_\Gamma(x, y, z)={F}_\Gamma(x, y, z)
+\sum_{i=1}^ra_if_i(x,y,z),
\]
where $f_i$ span the ring of polynomials in $(x, y, z)$ which do not
vanish when $\p_x{F}_\Gamma=\p_y{F}_\Gamma=\p_z{F}_\Gamma=0$.
The dimension $r$ of this ring is equal to the number of non-trivial
conjugacy classes of $\Gamma$ which is $k-1, k+1, 6, 7$ and $8$
respectively \cite{B70}.
Kronheimer \cite{K1,K2} proved that for each $\Gamma$ a unique hyper-K\"ahler metric
exists on a minimal resolution ${\cal M}$, and that this metric
is precisely the ALE metric with $\R^4/\Gamma$ as its infinity.
His construction was a combination of the hyper-K\"ahler quotient
\cite{HKLR87} with twistor theory.

In each case the twistor space is the three dimensional hyper-surface
$\widetilde{F}_\Gamma(x, y, z, \ll)=0$
in the rank-three bundle
$\OO(p)\oplus\OO(q)\oplus\OO(r)\rightarrow \CP^1$. Now
$
x(\ll)\in \OO(p), y(\ll)\in Q(q), z(\ll)\in O(r)
$
are polynomials in $\ll$, $f_i=f_i(x,y,z)$, and $a_i=a_i(\ll)$.
Therefore 
\[
{\cal PT}\longrightarrow\OO(p),\qquad {\cal
PT}\longrightarrow\OO(q),\qquad
{\cal PT}\longrightarrow\OO(r),
\]
and Lemma \ref{kstf} implies that  
the corresponding hyper-K\"ahler metrics admits three commuting hidden
symmetries, and 
the heavenly equation (\ref{secondeq}) reduces to an ODE. 

The degrees $p,q$ and $r$ are such that $\widetilde{F}_\Gamma(x, y, z,
\ll)$ is a function homogeneous of some degree $s$.
Therefore
\[
\widetilde{F}_\Gamma:\OO(p)\oplus\OO(q)\oplus\OO(r)\rightarrow \OO(s).
\]
To determine the integers $p, q, r, s$ 
take the determinants of the above, and notice that  
the normal bundle to an $\OO(1)\oplus\OO(1)$ section of ${\cal 
PT}\longrightarrow\CP^1$ will have the Chern class
$p+q+r-s=2$. This gives us the following
\begin{eqnarray}
\label{tabela}
A_k & & {\cal PT}=\{(x, y, z, \ll)\in \OO(k)\oplus\OO(k)\oplus\OO(2)
\longrightarrow \CP^1,\nonumber\\
&& xy-z^k-a_1z^{k-2}-...a_{k-1}=0\},\nonumber\\
D_{k-1} & & {\cal PT}=\{(x, y, z, \ll)\in \OO(2k)\oplus\OO(2k-2)\oplus\OO(4)\longrightarrow \CP^1,\nonumber\\
&&x^2+y^2z+z^{k}+a_1y^2+a_2y+a_3z^{k-2}+...+a_{k}z+a_{k+1} =0\},\\
E_6&&{\cal PT}=\{(x, y, z, \ll)\in \OO(12)\oplus\OO(8)\oplus\OO(6)
\longrightarrow \CP^1,\nonumber\\
&&x^2+y^3+z^4
+y(a_1z^2+a_2z+a_3)+a_4z^2+a_5z+a_6=0\}\nonumber\\
E_7&& {\cal PT}=\{(x, y, z, \ll)\in\OO(18)\oplus\OO(12)\oplus\OO(8)\longrightarrow \CP^1,\nonumber\\
&&x^2+y^3+yz^3+
y^2(a_1z+a_2)+y(a_3z+a_4)+a_5z^2+a_6z+a_7=0\}\nonumber\\
E_8&&{\cal PT}=\{(x, y, z, \ll)\in\OO(30)\oplus\OO(20)\oplus\OO(12)
\longrightarrow \CP^1,\nonumber\\
&&x^2+y^3+z^5+y(a_1z^3+a_2z^2+a_3z+a_4)+a_5z^3+a_6z^2+a_7z+a_8=0\}\nonumber
\end{eqnarray}
Note that these spaces are not quite the full non-singular twistor
space as there will be singular points where 
$\widetilde{F}_\Gamma$ vanishes together 
with its first derivatives.  These singularities can, however, be
resolved, see \cite{K2}.

We observe that these twistor spaces have projections onto $\OO(2n)$ for
$2n=p,q,r$ and this corresponds to the existence of three independent
commuting tri-holomorphic hidden symmetries.  
The simplest description along the
lines of \S\ref{patching_functions} arises for the lowest value of
$n$, i.e., when we project onto the $z$ coordinates in the above
construction.  It is clear from the above formulae that the fibres of
this projection are affine conics in the $A_k$ and $D_k$ cases, and
affine elliptic curves in the $E_k$ cases.

From now on we shall drop the subscript $\Gamma$, because the construction
we shall describe applies to all cases.
These twistor spaces can all be described as affine line bundles over
$\OO(r)$ in effect by uniformizing  the affine conics or elliptic curves
that make up the fibres of $\cPT\rightarrow \OO(r)$.  The
affine line bundle is determined by a cohomology class in
$H^1(\OO(r),\OO(2-r)$ which can be evaluated as a linear field on
$\C^{r+1}$.  The ALE space can then be realised as (a branched cover of) 
the zero set of this linear field in the real slice $\R^{r+1}$ as in
theorem (\ref{maintheorem}).

To make the description more concrete, we now find a patching
description of the relevant cohomology class in $H^1(\OO(r),\OO(2-r))$.
We exclude the curve(s) on which both $\widetilde{F}_x= 0$ and
$\widetilde{F}_y= 0$ so that we can cover the twistor space by the two
open sets $U, \widetilde{U}$ such that $\widetilde{F}_x\neq 0$ in $U$
and $\widetilde{F}_y\neq 0$ in $\widetilde{U}$.  We use $(y, z, \ll)$
and $(x, z, \ll)$ as local coordinates in $U$ and $\widetilde{U}$
respectively.  The symplectic form $\Sm$ on each fibre of ${\cal
  PT}\rightarrow\CP^1$ is given by
\[
\frac{\d y\wedge\d z}{\widetilde{F}_x}\;\;\mbox{in}\;U
,\qquad\mbox{or} 
\qquad\frac{-\d x\wedge\d z}{\widetilde{F}_y}\;\;
\mbox{in}\;\widetilde{U}.
\]
These arise from the formula 
$\Sigma=\oint \d x\wedge \d y \wedge \d z/\widetilde{F}$ with the contour being a small circle surrounding 
$\widetilde{F}=0$.  The
global homogeneous function $z$ gives rise to a homogeneity $r-2$
Hamiltonian vector field $X_z$ tangent to the fibres of ${\cal
  PT}\rightarrow\OO(r)$.  From the formula $X_z\hook\Sm=\d z$ we
deduce that
\[
X_z=\widetilde{F}_x\frac{\p}{\p y}
\;\;\mbox{in}\;U,\qquad
X_z=-\widetilde{F}_y\frac{\p}{\p x}
\;\;\mbox{in}\;\widetilde{U}.
\]
We now introduce new coordinates $(u, z, \ll)$ and $(\tilde{u}, z, \ll)$
on $U$ and $\widetilde{U}$ respectively, where the fibre coordinates
in $u$ in $U\rightarrow\OO(r)$ and $\tilde{u}$ in
$\widetilde{U}\rightarrow\OO(r)$ satisfy $X_z(u)=X_z(\tilde{u})=1$.
Therefore
\[
u(y, z)=\int_{\sigma}^{\widetilde{F}_x=1} \frac{\d y}{\widetilde{F}_x},
\qquad
\tilde{u}(x, z)=-\int_{\sigma}^{\widetilde{F}_y=1}\frac{\d x}{\widetilde{F}_y} 
\]
for some $\sigma$, 
and the patching function is given on the overlap by
\be
\label{patch_fun}
f(z,\ll)=u-\tilde{u}=\int_{\widetilde{F}_y=1}^{\widetilde{F}_x=1} \frac{\d y}{\widetilde{F}_x}.
\ee
In the above formula $x$ should be determined in terms of $(y, z, \ll)$
using $\widetilde{F}=0$ before the integral is evaluated. The upper
and lower limits will then involve $y=y(z, \ll)$. 

In the case
of $A_k$ ALE space we can assume 
that\footnote{See \cite{H79} and \cite{B70} for further discussion of
  this point.}
\[
z^k+a_1z^{k-2}+...+a_{k-1}=\prod_{j=1}^k(z-p_j(\ll)),\qquad\mbox{where}\;\;
p_j\in\Gamma(\OO(2)).
\]
A simple integration yields $f=\ln{\prod_{j=1}^k(z-p_j(\ll))}$, and
from Proposition \ref{two_methods} we find 
\[
G=\sum_{j=1}^k(z-p_j)(\ln{(z-p_j)-1)}.
\]
For $D_{k-1}$ we can redefine $z$ and $a_j$ to get rid of terms
linear and quadratic in $y$, and write
\[
\widetilde{F}=x^2-y^2z-\prod_{j=1}^k(z-q_j(\ll))\qquad\mbox{where}\;\;
q_j\in\Gamma(\OO(4)).
\]
Now (\ref{patch_fun}) yields
\[
f=\frac{1}{\sqrt{z}}\ln{\frac{\big(z-4z\prod_{j=1}^k(z-q_j(\ll))\big)^{1/2}+\sqrt{z}}{\big(1+4z\prod_{j=1}^k(z-q_j(\ll))\big)^{1/2}-1}}.  
\]
In the remaining cases $E_6, E_7$ and $E_8$ the fibres of 
${\cal PT}\longrightarrow \OO(r)$ are elliptic curves
$x^2=4y^2+g_1y+g_2$ (in case of $E_7$ one needs to redefine $y, z, a_i$ to
obtain this canonical form). The periods $g_1, g_2$ are polynomials in $z$
of order less or equal to $5$ which can be determined form (\ref{tabela}).
The fibres can therefore be parametrised by the Weierstrass elliptic function.
The cohomology class is represented by an elliptic integral
\[
f=\frac{1}{2}\int_{y_0}^{y_1}\frac{\d y}{\sqrt{4y^3+g_1y+g_2}},
\]
where $y_1$ and $y_0$ are roots of $4y^3+g_1y+g_2-1/4=0$ and $12y^2+g_1-1=0$
respectively.

One can now, in principle, take these cohomology classes and integrate
them to obtain a linear field with $r-3$ components on $\C^{r+1}$ the
vanishing of which will determine the complexified ALE space as a
submanifold.  This above description is not completely satisfactory
for two related reasons.  Firstly the description of the ALE space
will not be global; the projection from the true ALE space to
$\C^{r+1}$ can be many to one, and can have irregular values.  
Secondly, the limits of integration above defining the cohomology
classes actually branch and are not completely well defined. 

Further work is required to make this a useful description of ALE spaces.
It seems likely that these are the only complete hyper-K\"ahler metrics
with three tri-holomorphic hidden symmetries.

\section{Hierarchies for the generalised conformal anti-self-duality
  equations}\label{gcasdhier}
In this section we extend the concept of a hierarchy from that of
\cite{DM00} for the 4-dimensional hyper-K\"Ahler equations to a
generalisation of the conformal anti-self-duality equations (and in the
process give new and more geometric formulations for the
hyper-K\"ahler hierarchy than in \cite{DM00}).  The guiding motivation
for these definitions come from the twistor theory.  However, we first
define the various concepts in space-time terms, and then discuss the
twistor theory subsequently.  We shall, for convenience,
work in the holomorphic category.  Real versions of the various
structures and equations can then be obtained subsequently by
demanding the existence of an anti-holomorphic involution $\sigma$
fixing a real slice and with specified action on the various
geometric structures.

We will abbreviate the term conformal anti-self-duality to CASD and
generalised CASD to GCASD.  Unfortunately this terminology is non
standard but is designed to be consistent with the corresponding
discussion for the anti-self-dual Yang-Mills equations given in
\cite{MW96}.  The generalisation of the CASD case is a mild
generalisation of quaternionic structures discussed in \cite{Sa82} and
have been termed paraconformal structures, see \cite{BE} and Grassman
structures, \cite{A,BC} where many properties, including the twistor
theory, of these spaces are studied.  Here we shall refer to them as
generalised CASD, GCASD, spaces.  The hierarchies defined here are a
special case of the ${\mathcal P}$-structures of Gindikin, \cite{G90}
and references therein.

\begin{defi}
A solution to the GCASD hierarchy consists of the data
$(\cM,\bbS,\tbbS,e^{AA'_1\cdots A'_n})$ defined as follows: $\cM$ is a
manifold of dimension $r(n+1)$, $\bbS$ and $\tbbS$ are vector bundles
of rank $r$ and $2$ respectively, we use abstract indices $A$ and $A'$
to denote membership of $\bbS$ and $\tbbS$ respectively; when realised
concretely, $A=0,1, \cdots , r-1$ and $A'=0',1'$.  The indexed 1-form
$e^{AA'_1\cdots A'_n}$, symmetric over its primed indices, determines an
isomorphism $T\cM=\bbS\otimes \odot^n \tbbS$ at every point.

An element $\pi_{A'}$ of $\tbbS^*$ at $m\in \cM$ determines an $rn$--plane
element 
$$
z(m)_\pi=\{V\in T_m\cM, V\hook e^{AA'_1\cdots
A'_n}\pi_{A'_1}\cdots \pi_{A'_n}=0\}.
$$
Such an $rn$-plane element will be said to be an $\alpha$-plane
element at $m$.  An $\alpha$-surface is an $rn$-dimensional surface
whose tangent space defines an $\alpha$-plane element at each of its
points.

The GCASD hierarchy equations are the requirement that there exists a
full family of $\alpha$-surfaces, with a unique $\alpha$-surface
through each $z(m)_\pi$.
\end{defi}

The notation derives from the identification of these bundles with the
spin bundles of a conformal structure in 4-dimensions, $r=2$, $n=1$.
It will also be convenient to introduce a `clumped' index $i$ for the
$n+1$-dimensional vector space $\odot^n\tbbS$.  When the indices are
realised concretely by a choice of a frame for $\tbbS$ with components
labelled by $0$ and $1$, there is a standard correspondence between
the $i$th component for the clumped index, and the component with $i$
$1$s and $n-i$ $0$s, so that $i$ naturally goes from $0$ to $n$.

We now assume that we have a solution to the CASD hierarchy so that we
have a full complement of $\alpha$-surfaces and that, shrinking $\cM$ to
a convex neighbourhood of a point if necessary, the space of these
$\alpha$-surfaces is a manifold.  We can then define

\begin{defi}
The space of such $\alpha$-surfaces will be called the twistor space
and is denoted $\PT$.
\end{defi}

Twistor space is an $r+1$ dimensional complex manifold.  
\begin{theo}
The twistor space determines and is determined by the GCASD
hierarchy.  The correspondence is stable under small deformations of
the complex structure of $\PT$ or of the GCASD hierarchy.
\end{theo}

\Proof
This is a straightforward extension of Penrose's nonlinear graviton
construction.  The correspondence can be studied by means of the double
fibration 
\begin{eqnarray*}
&P\tbbS&\\
p\swarrow&&\searrow q\\
\cM\quad && \quad \PT\, .
\end{eqnarray*}
Points $m\in M$ correspond to rational curves $L_m:=q(p^{-1}(m))
\equiv \CP^1$ in $\PT$.  The normal bundle of these rational curves is
$N=\bbS \otimes \OO(n)$ where $\OO(n)$ is the line bundle of Chern class
$n$ on $\CP^1$;  this follows from the fact that, since the tangent
space of $z(m)_\pi$ is the kernel of
$e^{AA'_1\cdots A'_n}\pi_{A'_1}\cdots \pi_{A'_n}$, the section of the
normal bundle corresponfing to a vector $V$ can be identified with
$V\hook e^{AA'_1\cdots A'_n}\pi_{A'_1}\cdots \pi_{A'_n}$, a function
with values in $\bbS$ with homogeneity $n$.  However, sections of
$\OO(n)$ can be identified with functions homogeneous degree $n$, and
so the normal bundle is $\bbS\otimes\OO(n)$ as claimed.

With knowledge of the normal bundle, Kodaira theory can now be applied
and shows that, since $H^1(\CP^1,N)=H^1(\CP^1,End(N))=0$, the moduli
space of curves has $dim (H^0(\CP^1,\bbS\otimes\OO(n))=
r(n+1)$ dimensions, contains $\cM$ and $T_m\cM\equiv
\Gamma(\bbS\otimes\OO(n))\equiv \bbS\otimes\odot^n\tbbS$. Points of
$\PT$ clearly then correspond to integrable $\alpha$-surfaces in $\cM$.
Kodaira theory also provides the stability of the correspondence under
small deformations.  See \cite{ME1,ME2} for general constructions that
apply to these situations.  \koniec

\noindent
{\bf Remark:} The ${\mathcal P}$-structures of Gindikin are more
general but can be understood easily in this context as arising
naturally on moduli spaces of rational curves in some complex manifold
whose normal bundles are $\OO(k_1)\oplus \cdots \OO(k_r)$ with the $k_i$
not being required to be equal.  Such prescriptions for the normal
bundle are unstable under deformations of the underlying complex
manifold unless no two of the $k_i$ differ by more than one.  Under
deformations of the complex structure, in the moduli space of such
rational curves, the normal bundle will jump so that a dense open set
will be the stable case where the normal bundle will be
$E\otimes\OO(k)\oplus F\otimes \OO(k+1)$ with $E$ and $F$ trivial
bundles.  

This additional generality can be important.  For example, if we wish
to discuss generalisations of the Ward construction, the twistor space
has the structure of a holomorphic vector bundle over a lower
dimensional space.  In this case, the normal bundle of the rational
curves along the fibres is usually taken to have degree $0$, whereas
the normal bundle of its projection into the bases will usually be
required to have a higher degree normal bundle.

\bigskip

These geometric structures fall into the category of involutive
$G$-structures studied by Merkulov \cite{ME1,ME2}.  In particular, one
can exploit his theorems to deduce the existence of connections
compatible with the geometric structure.  However, for the most part,
they must have torsion, although they fall into Merkulov's category of
`$G$-structures with very little torsion', \cite{ME2}.

\begin{lemma}\label{conn}
For $r\geq 2$, $n\geq 1$, there exists connections on $\bbS$ and
$\tbbS$ such that the induced connection on $TM$ has torsion with
non-vanishing irreducible parts only in $\bbS\otimes
\odot^2\bbS^*\otimes\odot^{n-2}\tbbS$ and
$\bbS^*\otimes\odot^{n-4}\tbbS$ where we take $\odot^n \tbbS=\C$ for
$n=0$ or zero for $n<0$ .  There exists a unique choice for such a
connections when $n>1$ and unique up to a 1-form for $n=1$ (which can
be taken to be exact with appropriate choices). 
\end{lemma}

\Proof Merkulov reformulates the moduli spaces considered above as
Legendrian moduli spaces of holomorphically embedded
$\CP^1\times\CP^{r-1}$s in the projective cotangent bundle $PT^*\PT$
of twistor space.  A rational curve, $\CP^1$ in $\PT$ determines its
projective conormal bundle in $PT^*\PT$, i.e., the 1-forms up to scale
that annihilate the tangent space of $\CP^1$.  This correspondence is
studied by means of the following double fibration
\begin{eqnarray*}
&P\tbbS\times P\bbS &\subset P(T^*M)\\
\mu \swarrow&&\searrow \nu\\
M\quad && \quad P(T^*\PT)\, .
\end{eqnarray*}
Merkulov's method uses the contact line bundle $L$ which is the dual
to the tautological line bundle $T^*\PT\rightarrow P(T^*\PT)$.  On
restriction to a $\CP^1\times \CP^{r-1}$, it gives $\OO(n,1)$ where
$\OO(p,q)$ is the product of the pullback of $\OO(p)$ from $\CP^1$ with
the pullback of $\OO(q)$ frpm $\CP^{r-1}$.

Merkulov shows that the minimal torsion of an affine connection
preserving the $G$-structure (or obstruction to obtaining a torsion
free connection preserving the tensor decomposition of the tangent
space) is then measured by a geometrically obtained class in
$H^1(\CP^1\times \CP^{r-1},L\otimes\odot^2(J^1L)^*)$ and the freedom
in the resulting connection is given by $H^0(\CP^1\times
\CP^{r-1},L\otimes\odot^2(J^1L)^*)$.

These groups can be computed as follows.  The first jet of a section
of $\OO(n)$, $n\neq 0$, at a point of $\CP^{r-1}$ can be encoded into
the derivative of a homogeneous degree $n$ function with respect to
the $r+1$ homogeneous coordinates.  The value of the function is then
retrieved from this by Euler's homogeneity equations, $\pi_A\p f/\p
\pi_A=n f$.  Thus, the sheaf $J^1L$ on $\CP^1\times \CP^r$ can be
understood as the kernel
$$
0\rightarrow J^1L\longrightarrow \bbS^*(n-1,1)\oplus \tbbS^*(n,0)
\stackrel{(\pi_A, -\pi_{A'})}{\longrightarrow} \OO(n,1) \rightarrow 0
$$
since in the third map, we are imposing the requirement that the Euler
homogeneity relation for each factor leads to the same value for $f$.

The cohomoolgy groups $H^i(\CP^1\times
\CP^{r-1},L\otimes\odot^2(J^1L)^*)$ can therefore be calculated by
consideration of the long exact cohomology sequence arising from the
short exact sequence
\begin{eqnarray*}
0\longrightarrow
\bbS^*(-n,0)\oplus\tbbS^*(1-n,-1)
\stackrel{\odot(\pi_B,\pi_{B'})}{\longrightarrow} 
\odot^2\bbS^*(-n,1)\oplus \bbS\otimes\tbbS(1-n,0)\oplus  
\odot^2\tbbS(2-n,-1) \\
\hfill\rightarrow L\otimes\odot^2(J^1L^*)\rightarrow 0
\end{eqnarray*}
where $\pi_A$ and $\pi_{A'}$ are the homogenous coordinates on
$\CP^{r-1}$ and $\CP^1$ respectively.  
\koniec

Note that the Merkulov framework is not quite equivalent to ours in
the sense that it only requires knowledge of the total space $P(T^*\PT)$ but
does not require that it be realised as the projective cotangent
bundle of some $\PT$.  The results are only inequivalent for $r=1$,
$n<3$ and $r=2$, $n=1$.  In these cases the results are well known,
for example, the full theory of the latter case goes back to
Penrose 1976, \cite{Pe76}.  The Merkulov framework doesnt see the
curvature conditions that arise from existence of $\PT$, but gives the
correct result for the existence of and freedom in choosing compatible
torsion-free affine connections.

The cases $r=1$ are also well known, but for $n=1,2$ do not fall
satisfactorily into the Merkulov framework.  For $n=1$, there is a
projective structure, i.e., an equivalence class of torsion-free
connections that share the same unparametrised geodesics, with freedom
given by a 1-form.  For $n=2$ there exists a unique torsion-free
connection compatible with a conformal structure.  For $n=3$ the
connection is still torsion-free, but not subsequently for higher $n$.

The general $n=1$ case was studied in \cite{Sa82,BE}.  It also follows
from the calculations of \cite{BE} that the torsion must be non zero
for a non flat structure in the $n>1$, $r>2$ cases as a consequence of
the fact that in these cases the decomposition of the tangent space as
a tensor product of $\bbS$ with $\odot^n\tbbS$ determines a
paraconformal structure in which both factors have dimension greater
than two, and in that case the torsion-free condition implies
flatness.

\begin{lemma}
For $r>1$, the requirement of uniqueness for the $\alpha$ surface
through $z(m)_\pi$ is redundant.
\end{lemma}

\Proof
The integrability equations in particular give a propagation equation
for $\pi$ across the $\alpha$-surface. \koniec

\smallskip

There are a number of specializations of the GCASD equations:
hypercomplex, scalar-flat K\"ahler, Einstein, hyper-K\"ahler.  The
hypercomplex and hyper-K\"ahler cases have straightforward extensions to
the hierarchy.

\begin{itemize}
\item The hyper-complex case for $r$ even, $n=1$, where there exists a
flat connection on $\tbbS$ such that the distribution $D$ on $P\tbbS$
is horizontal.  This is equivalent to the existence of a fibration
$\PT\mapsto \CP^1$.  This condition (flatness of the induced
connection of \ref{conn} or a fibration of the associated twistor
space over $\CP^1$) can clearly be imposed consistently on any
$\cM_{r,n}$ to give a hypercomplex hierarchy.





\item In case of the hyper-K\"ahler hierarchy, we require that there
exists a connection compatible with \ref{conn} that induces a flat
connection on $\tbbS$ and preserves skew forms
$\varepsilon_{AB}$ on $\bbS$ and $\varepsilon_{A'B'}$ on $\tbbS$ such
that the forms $\varepsilon_{AB}e^{A(A'_1\cdots A'_n}\wedge
e^{B'_1\cdots B'_n)B}$ are closed.  This then implies that
$\eta=\varepsilon^{A'B'}\pi_{A'}D\pi_{B'}$ and $\eta\wedge
\varepsilon_{AB}\pi_{A'_1}\cdots \pi_{A'_n}e^{AA'_1\cdots A'_n}\wedge
\pi_{B'_1}\cdots \pi_{B'_n}e^{BB'_1\cdots B'_n}$ descend to $\PT$, in
such a way that $\eta$ is the annihilator of an integrable
distribution determining a fibration over $\CP^1$.

\end{itemize}
\subsection{Reality structures}
The imposition of reality conditions is standard; it is imposed by
requiring the existence of an anti-holomorphic involution $\sigma$ on
$\cM$ that fixes a real slice and preserves the geometric structures
(i.e., sends $\alpha$-surfaces to $\alpha$-surfaces).  In the
hyper-K\"ahler case, we can talk in terms of the signature of the
associated metric (although the following conditions can be applied
more generally).  For Euclidean signature, we require that it induces
a quaternionic involution on $\tbbS$ given by $\sigma^2=-1$.
In particular there are no non-zero fixed points.  It will then also induce a
quaternionic involution on $\bbS$ which will have to be even
dimensional and we must also require that the Hermitian form
$\pi^A\hat \pi^B\varepsilon_{AB}$ be definite (it is
trivially definite for $r=2$).  For non-Euclidean signature we can
have different signatures for $\pi^A\hat \pi^B\varepsilon_{AB}$, or
impose a conjugation whose action on $\bbS$ snf $\tbbS$ is an ordinary
complex conjugation.

The conjugation will also lead to an anti-holomorphic involution on
the twistor space, without fixed points in the quaternionic case, and
with a fixed real slice otherwise.  Points of the real slice of $\cM$
will then correspond to $\sigma$ invariant rational curves.

\subsection{Embedding into hierarchies}
In the usual definition of a hierarchy, the hierarchy is an
overdermined, but compatible system of equations which contains the
original system.  A given solution to the original system may not
actually extend to a solution of the hierarchy in general (there can
be obstructions, see p249 and p253 footnote 2 of \cite{MW96} for some
discussion of this behaviour for the Drinfeld Sokolov hierarchies).
Furthermore, if such an extension does exist, it will not in general
be unique without the imposition of boundary conditions.

Our definition of a hidden symmetry in  Section 4  
require the existence of an extension of a solution to the original
equation to the hierarchy that happens to admit a symmetry, but only
when thought of as a solution to the hierarchy.  

We state the embedding definition in slightly greater generality as
for one hierarchy into another:

\begin{defi}
A solution $\cM_{r_1,n_1}$ to the GCASD hierarchy embeds into another
solution $\cM_{r_2,n_2}$, $n_1<n_2$, $r_1\leq r_2$ if $\cM_{r,n_1}$ embeds
into $\cM_{r,n_2}$ as a manifold in such a way that the
$\alpha$-surfaces of $\cM_{r,n_2}$ intersect $\cM_{r,n_1}$ in the
$\alpha$-surfaces of $\cM_{r,n_1}$ and all $\alpha$-surfaces of
$\cM_{r,n_1}$ arise in this way.
\end{defi}
Due to a remarkable theorem of Bernstein \&
Gindikin, \cite{G90}, the twistor characterisation of such an
embedding in the most interesting case, $r_1=r_2$, is remarkably simple:

\begin{theo}[Bernstein \& Gindikin]
A solution $\cM_{r,n_1}$ to the GCASD hierarchy embeds into $\cM_{r,n_2}$
iff the twistor space $\PT_{n_1}$ for $\cM_{r,n_1}$ is obtained from
that, $\PT_{n_2}$ for $\cM_{r,n_2}$ by choosing submanifolds $\Gamma_1$,
$\Gamma_2$, $\cdots$ of codimension greater than one and $S_1$, $S_2$,
$\cdots$ of codimension$=1$ and blowing up along each $\Gamma_1$,
$\Gamma_2$, $\cdots$ and taking a branched covers branching with some
multiplicity over each $S_i$.
\end{theo}

In \cite{DM00} (see also Section 3), 
we embedded $\cM_{2,n_1}$ into $\cM_{2,n_2}$ by blowing
up the twistor space $\PT_{n_2}$ at one point $n_2-n_1$ times.

Note that it is natural in the quaternionic cases, or in the
hyper-CASD cases to require that the additional structures be
compatible.  This is straightforward in the hypercomplex case in which
one wishes the embedded twistor space to inherit a fibration over
$\CP^1$, but when line bundle valued forms need to be pulled back
also, there is the problem that the forms that have been pulled back
will in general take values in an inappropriate line bundle unless
particular care has been taken.

\subsection{Symmetries and hidden symmetries}
We can define a symmetry for a GCASD structure to be the requirement
that the twistor space admits a global holomorphic vector field $\cK$.
This will in turn determine global holomorphic vector fields
$\tilde{K}$ on the correspondence space and $K$ on $\cM$, such that
$\tilde{K}$ projects to $K$.  The essential requirement on $\tilde{K}$
will be that it preserves the twistor distribution and this will lead
to a generalisation of the conformal Killing vector equations on $K$
whose precise form will depend on $r$ and $n$.  

Clearly the concept of hidden symmetry can be applied as before but
with greater generality; a solution $\cM_{r,n}$ admits a hidden
symmetry if it can be embedded as above into an $\cM_{r,m}$ that
admits an explicit symmetry.  This will, as in the proof of
proposition \ref{hiddensym}, lead to a global vector field on
$\cPT_{r,n}$ with values in a line bundle $\Lie$ of degree $m-n$ (the
restriction of the canonical bundle of $\cPT_{r,m}$ tensored with the
inverse of that of $\cPT_{r,n}$). There will also be a generalisation
of theorem \ref{maintheorem}: this global vector field with values in
$\Lie$ will lead to the realisation of $\cPT_{r,n}$ as the total space
of an affine line bundle, with underlying translation bundle $\Lie^*$,
over some reduced twistor space which will generically be
$\cPT_{r-1,l}$ where $l=n+m/(r-1)$, if $l$ is an integer, although if
$l$ is fractional, or in non-generic situations, the normal bundle of
lines in the reduced twistor space must be $\OO(k_1)\oplus \cdots
\oplus \OO(k_{r-1})$ with $\sum k_i = (r-1)n +m$.  Thus the original
$\cPT_{r,n}$ with a hidden symmetry can be determined in terms of a
lower dimensional twistor space together with a linear cohomology
class on that space.

\section*{Acknowledgements}
We are grateful to Nigel Hitchin for helpful comments and for
showing us his notes on $D_k$ ALE manifolds. 
We would also like to thank Toby Bailey, Roger Bielawski, David Calderbank,
Simon Gindikin and George Sparling for useful discussions.
Both authors are members of EDGE, Research Training Network 
HPRN-CT-2000-00101. MD was partly supported by NATO grant PST.CLG.978984.
\section*{Appendix}
In four complex dimensions orthogonal transformations decompose into
products of ASD and SD rotations 
$ \label{basicisom}
SO(4, \C)=(SL(2, \C)\times \widetilde{SL}(2, \C))/\Z_2. $ 
The spinor calculus in four dimensions is based on this isomorphism.  We
use the conventions of Penrose and Rindler \cite{PR86}.  
 The tangent space at each point of
${\cal M}$ is isomorphic to a tensor product of the two spin spaces 
$
\label{isomorphism}
T^a{\cal M}=S^A\otimes S^{A'}.
$
Spin dyads $(o^A, \iota^{A})$ and $(o^{A'}, \iota^{A'})$ span $S^A$
and $S^{A'}$ respectively.  The spin spaces $S^A$ and $S^{A'}$ are
equipped with symplectic forms $\varepsilon_{AB}$ and
$\varepsilon_{A'B'}$ such that
$\varepsilon_{01}=\varepsilon_{0'1'}=1$.  These anti-symmetric objects
are used to raise and lower the spinor indices via 
$\iota_A=\iota^B\varepsilon_{BA}, \iota^B=\varepsilon^{AB}\iota_B$
Let $\Gamma_{AB}$ and $\Gamma_{A'B'}$ are the  $SL(2, \C)$
and $\widetilde{SL}(2, \C)$ spin connections.
The curvature of the unprimed spin connection
$
{R^A}_B=\d{\Gamma^A}_B+{\Gamma^A}_C\wedge{\Gamma^C}_B
$
decomposes as
\[
{R^A}_B={C^A}_{BCD}\Sm^{CD}+(1/12)R{\Sm^{A}}_{B}+{\Phi^A}_{BC'D'}\Sm^{C'D'},
\]
and similarly  for ${R^{A'}}_{B'}$. Here $R$ is the Ricci scalar, 
$\Phi_{ABA'B'}$ is the trace-free part of the Ricci tensor $R_{ab}$,
and $C_{ABCD}$ is the ASD part of the Weyl tensor
\[
C_{abcd}=\varepsilon_{A'B'}\varepsilon_{C'D'}C_{ABCD}+
\varepsilon_{AB}\varepsilon_{CD}C_{A'B'C'D'},
\]
and the  two forms $\Sm^{A'B'}$ span the three dimensional space of SD
two forms.

Given a complex four-dimensional manifold ${\cal M}$ with curved
metric $g$, a twistor in ${\cal M}$ is an $\a$-surface, i.e.\ a null 
two-dimensional surface whose tangent space at each point is an $\a$
plane(a  null 
two-dimensional plane with a SD bi-vector).  
There are Frobenius integrability conditions for the existence
of such $\a$-surfaces through each $\a$-plane element at each point
and these are equivalent, after some calculation, to the vanishing of
the self-dual part of the Weyl curvature, $C_{A'B'C'D'}$.  Thus, given
$C_{A'B'C'D'}=0$, we can define a twistor space ${\cal PT}$ to be the
three complex dimensional manifold of $\a$-surfaces in ${\cal M}$.  If
$g$ is also Ricci flat then ${\cal PT}$ has further structures which
are listed in the Nonlinear Graviton Theorem:
\begin{theo}[Penrose \cite{Pe76}]
\label{Penrose}
There is a 1-1 correspondence between complex ASD vacuum metrics
on complex four-manifolds 
and three dimensional complex manifolds ${\cal PT}$ such that
\begin{itemize}
\item There exists a holomorphic projection $\mu:{\cal
PT}\longrightarrow \CP^1$
\item ${\cal PT}$ is equipped with a four complex parameter
family of sections of $\mu$ each with a normal bundle $
{\cal O}(1)\oplus{\cal O}(1)$, (this will follow from the existence of
one such curve by Kodaira theory),
\item Each fibre of $\mu$ has a symplectic  structure
$
\Sm_{\ll} \in \Gamma(\Lambda^2(\mu^{-1}(\lambda))\otimes {\cal O}(2)),
$
where $\ll\in \CP^1$.
\end{itemize}
\end{theo}
To obtain real metrics on a real four-manifold, we can require further that the
twistor space admit an anti-holomorphic involution.

The correspondence space ${\cal F}={\cal M}\times\CP^1$ is
coordinatized by $(x,\lambda)$, where $x$ denotes the coordinates on
$\cal M$ and $\ll$ is the coordinate on $\CP^1$ that parametrises the
$\a$-surfaces through $x$ in $\cal M$.  We represent $\cal F$ as the
quotient of the primed-spin bundle $S^{A'}$ with fibre coordinates
$\pi_{A'}$ by the Euler vector field $\Upsilon=\pi^{A'}/\p
\pi^{A'}$. We relate the fibre coordinates to $\lambda$ by
$\lambda=\pi_{0'}/\pi_{1'}$.
A form with values in the line bundle ${\cal O}(n)$ on $\cal F$ can be
represented by a homogeneous form $\alpha$ on the non-projective spin
bundle satisfying $
\Upsilon\hook\alpha=0\, , \qquad {\cal L}_{\Upsilon} \alpha=n\alpha .
$

The correspondence space has the alternate definition
\[
{\cal F}={\cal PT}\times {\cal M}|_{Z\in l_x}= {\cal M}\times\CP^1
\] 
where $l_x$ is the line in $\cal PT$ that corresponds to $x\in {\cal
M}$ and $Z\in\cal PT$ lies on $l_x$.  This leads to a double fibration
\be
\label{doublefib}
{\cal M}\stackrel{p}\longleftarrow 
{\cal F}\stackrel{q}\longrightarrow {\cal PT}.
\ee
Points in ${\cal M}$ correspond to
rational curves in ${\cal PT}$ with normal bundle ${\cal
  O}^{A}(1):={\cal O}(1)\oplus {\cal O}(1)$.  The normal bundle to
$l_x$ consists of vectors tangent to $x$ (horizontally lifted to
$T_{(x,\lambda)}{\cal F}$) modulo the twistor distribution. 
We have a sequence of sheaves over $\CP^1$
\[
0\longrightarrow D \longrightarrow \C^4 \longrightarrow
{\cal O}^A(1)\longrightarrow 0.
\]
The map $\C^4 \longrightarrow {\cal O}^A(1)$ is given by
$V^{AA'}\longrightarrow V^{AA'}\pi_{A'}$.  Its kernel consists of
vectors of the form $\pi^{A'}\lambda^A$ with $\lambda^A$ varying. The
twistor distribution is therefore $D=O(-1)\otimes S^{A}$ and so there
is a canonical $L_A\in\Gamma(D\otimes {\cal O}(1)\otimes S_{A})$
given by $L_A=\pi^{A'}\nabla_{AA'}$.
The projective twistor space $\cal PT$ arises as a quotient of ${\cal
F}$ by the twistor distribution. Functions on $\cal F$ which are
constant along $L_{A'}$ are called twistor functions.

\end{document}